\documentclass[10pt]{amsart}
\usepackage[margin=2.5cm]{geometry}
\usepackage{tikz-cd}
\usepackage{graphicx}
\usepackage{mathdots}		
\usepackage{amssymb}
\usepackage{amsmath}
\usepackage{relsize}
\usepackage{subfig, graphicx}
\usepackage{mathrsfs}
\usepackage{amssymb, amsthm, indentfirst}
\usepackage{relsize}
\usepackage{hyperref}
\usepackage{stmaryrd}
\usepackage{lineno}
\usepackage{mathabx}

\numberwithin{equation}{section}
\usepackage{color}
\theoremstyle{plain}
\newtheorem{thm}{Theorem}[section]

\newtheorem{lemma}[thm]{Lemma}
\newtheorem{prop}[thm]{Proposition}
\newtheorem{corollary}[thm]{Corollary}

\theoremstyle{definition}
\newtheorem{definition}[thm]{Definition}
\newtheorem*{example*}{Example}
\newtheorem{rmk}[thm]{Remark}
\def\Gal{\operatorname{Gal}}

\newtheorem*{hypothesis*}{Hypothesis}

\author{Shih-Yu Chen}
\address{Department of Mathematics, National Tsing Hua University, 101, Section 2, Kuang-Fu Road, Hsinchu 300, Taiwan, ROC}
\email{sychen.math@gmail.com}

\def\GL{{\rm{GL}}}
\def\GSp{{\rm GSp}}

\def\Sp{{\rm Sp}}

\def\A{{\mathbb A}}
\def\C{{\mathbb C}}
\def\E{{\mathbb E}}
\def\F{{\mathbb F}}

\def\R{{\mathbb R}}
\def\Q{{\mathbb Q}}
\def\Z{{\mathbb Z}}

\def\<{\langle}
\def\>{\rangle}

\def\x{\times}

\def\bp{\begin{pmatrix}}
\def\ep{\end{pmatrix}}

\def\<{\langle}
\def\>{\rangle}

\def\GL{\operatorname{GL}}

\def\GSp{\operatorname{GSp}}

\def\SO{\operatorname{SO}}
\def\Sp{\operatorname{Sp}}

\def\1{\mathbf{1}}

\def\itPi{\mathit{\Pi}}
\def\itPsi{\mathit{\Psi}}
\def\itSigma{\mathit{\Sigma}}
\def\itPhi{\Phi}

\makeatletter
\newcommand{\exterior}[1]{\mathop{\mathpalette\exterior@{#1}}}
\newcommand{\exterior@}[2]{%
  \raisebox{\depth}{%
  \fontsize{\sf@size}{0}%
  \m@th
  $\ifx#1\displaystyle\textstyle\else#1\fi\bigwedge$}%
  ^{\mspace{-2mu}#2}%
  \kern-\scriptspace
}
\makeatother

\title{Period relations between the Betti--Whittaker periods for $\GL_n$ under duality}

\begin{document}
\begin{abstract}
In this paper, under some regularity conditions, we prove a period relation between the Betti--Whittaker periods associated to a regular algebraic cuspidal automorphic representation of $\GL_n(\A)$ and its contragredient.
As a consequence, we prove the trivialness of the relative period associated to a regular algebraic cuspidal automorphic representation of $\GL_{2n}(\A)$ of orthogonal type, which implies the algebraicity of the ratios of successive critical $L$-values for ${\rm GSpin}_{2n}^* \times \GL_{n'}$ by the result of Harder and Raghuram.
\end{abstract}

\maketitle

\section{Introduction}

The purpose of this paper is to prove an automorphic analogue of a relation between the period invariants of motives under duality.
To be precise, let $M$ be a regular pure motive over $\Q$ of rank $n$ with coefficients in a number field $\E$.
The comparison isomorphism between the Betti and de Rham realizations of $M$ determines the period matrix $X_M$ which is an $n$ by $n$ matrix with coefficients in $\E\otimes_\Q\C$. 
In \cite{Yoshida2001}, Yoshida defined period invariants of $M$ as the evaluation of admissible polynomial functions at $X_M$. 
For instance, the Deligne's periods $\delta(M)$ and $c^\pm(M)$ introduced in \cite[\S\,2.6]{Deligne1979} are equal to $\det(X_M)$ and $f^\pm(X_M)$, where $f^+$ (resp.\,$f^-$) is the determinant of the upper left (resp.\,upper right) $d^+$ by $d^+$ (resp.\,$d^-$ by $d^-$) submatrix, and $d^\pm$ is the dimension of the $\pm$-eigenspace of the Betti realization of $M$ under the archimedean Frobenius action.
We consider a specific admissible polynomial function $f_{\rm BW}^\varepsilon$ for $\varepsilon \in \{\pm1\}$ if $n$ is even and $\varepsilon = d^+-d^-$ if $n$ is odd (cf.\,(\ref{E:motivic BW period}) below). 
For a regular pure motive $N$ of rank $n-1$ whose Hodge types are in good position relative to that of $M$, by the computation of Yoshida on Deligne's periods of tensor product motives (explicated by Bhagwat  \cite{Bhagwat2014} for regular motives), we have
\[
c^\pm(M\otimes N) = \delta(N)\cdot f_{\rm BW}^\varepsilon(X_M)\cdot f_{\rm BW}^{\varepsilon'}(X_N),\quad \varepsilon\varepsilon' = \pm1.
\]
In particular, the celebrated Deligne's conjecture \cite[Conjecture 2.8]{Deligne1979} then predicts that the algebraicity of critical values of the motivic $L$-function $L(M\otimes N,s)$ should be expressed in terms of $\delta(N)\cdot f_{\rm BW}^\varepsilon(X_M)\cdot f_{\rm BW}^{\varepsilon'}(X_N)$.
If we consider the dual motive $M^\vee$ of $M$, then Deligne observed in \cite[(5.1.7)]{Deligne1979} the period relation 
\[
c^\pm(M^\vee) = \delta(M)^{-1}\cdot c^\mp(M).
\]
Moreover generally, if $f$ is admissible of type $\{\underline{a};(k^+,k^-)\}$ and $f^\vee$ be its dual, 
then we have (cf.\,Lemma \ref{L:motivic} below)
\begin{align}\label{E:motivic period relation}
f^\vee(X_{M^\vee}) = \delta(M)^{-k^+-k^-}\cdot f(X_M).
\end{align}
In this paper we prove an automorphic analogue of the period relation (\ref{E:motivic period relation}) for $f=f_{\rm BW}^\varepsilon$.
More precisely, let $\itPi$ be a regular algebraic cuspidal automorphic representation of $\GL_n(\A)$, where $\A$ is the ring of rational adeles.
In \cite{Mahnkopf2005}, Mahnkopf introduced the Betti--Whittaker periods $p(\itPi,\varepsilon)$ for each admissible signature $\varepsilon$. For regular algebraic cuspidal automorphic representation $\itSigma$ of $\GL_{n-1}(\A)$ such that $(\itPi_\infty,\itSigma_\infty)$ is balanced, Raghuram proved in \cite{Raghuram2009} (generalizing the previous results of Harder \cite{Harder1983}, Hida \cite{Hida1994}, Kazhdan--Mazur--Schmidt \cite{KMS2000}, and Mohnkopf \cite{Mahnkopf2005}) that the algebraicity of critical values of the Rankin--Selberg $L$-function $L(s,\itPi \times \itSigma)$ can be expressed in terms of product of Betti--Whittaker periods $p(\itPi,\varepsilon)\cdot p(\itSigma,\varepsilon')$.
Note that the Rankin--Selberg $L$-function is equal to the (conjectural) motivic $L$-function $L(M_\itPi \otimes M_\itSigma,s)$ up to certain shift, where $M_\itPi$ and $M_\itSigma$ are the (conjectural) regular pure motives associated to $\itPi$ and $\itSigma$ proposed by Clozel \cite{Clozel1990}.
Therefore, we expect a relation between the Betti--Whittaker period $p(\itPi,\varepsilon)$ and the motivic period invariant $f_{\rm BW}^{\varepsilon}(X_{M_\itPi})$ (cf.\,\cite{HN2022}).
In particular, we anticipate a relation between $p(\itPi,\varepsilon)$ and $p(\itPi^\vee,\varepsilon)$ analogous to (\ref{E:motivic period relation}).
The aim of this paper is to establish this automorphic period relation Theorem \ref{T:Main 1} under some regularity conditions.
As a consequence of the period relation, we obtain Theorem \ref{T:Main 2} on the algebraicity of ratios of successive critical $L$-values for ${\rm GSpin}_{2n}^* \times \GL_{n'}$, where ${\rm GSpin}_{2n}^*$ is a quasi-split general spin group over $\Q$ of type $D_n$ or ${}^2D_n$.
This generalizes results of Bhagwat--Raghuram \cite{BR2020}.

\subsection{Main results}

Denote by $\A$ the ring of adeles of $\Q$.
Let $\itPi$ be a regular algebraic cuspidal automorphic representation of $\GL_n(\A)$. The regular algebraic condition implies that the finite part $\itPi_f$ is defined over its rationality field $\Q(\itPi)$ which is a number field. The archimedean component $\itPi_\infty$ determines a locally constant sheaf $\mathcal{M}$ of $\Q$-vector spaces on 
\[
\mathcal{S}_n = \GL_n(\Q)\backslash\left(\GL_n(\R)/\R_+{\rm SO}_n(\R)\times\GL_n(\A_f)\right)
\]
which is an inverse limit of locally symmetric spaces.
For each admissible signature $\varepsilon \in \{\pm1\}$, the representation $\varepsilon\times\itPi_f$ of $\pi_0(\GL_n(\R))\times\GL_n(\A_f)$ appears with multiplicity one in the cuspidal cohomology $H_{\rm cusp}^\bullet(\mathcal{S}_n,\mathcal{M}_\C)$ in bottom degree $\bullet = b_n = \lfloor \tfrac{n^2}{4} \rfloor$. We denote by $H_{\rm cusp}^{b_n}(\mathcal{S}_n,\mathcal{M}_\C)[\varepsilon\times\itPi_f]$ the corresponding isotypic component. 
Note that $\varepsilon$ can be arbitrary when $n$ is even, and is a uniquely determined signature $\varepsilon(\itPi_\infty)$ when $n$ is odd.
On the other hand, we have the Whittaker model $\mathcal{W}(\itPi_f)$ of $\itPi_f$ consisting of Whittaker functions of $\itPi_f$.
With respect to a choice of generator in the relative Lie algebra cohomology of $\itPi_\infty$ with suitable coefficients, the inverse map of taking global Whittaker functions of cusp forms in $\itPi$ induces a $\GL_n(\A_f)$-equivariant isomorphism
\[
\mathcal{W}(\itPi_f)\longrightarrow H_{\rm cusp}^{b_n}(\mathcal{S}_n,\mathcal{M}_\C)[\varepsilon\times\itPi_f].
\]
The (bottom degree) Betti--Whittaker period $p(\itPi,\varepsilon)\in\C^\times/\Q(\itPi)^\times$ of $\itPi$ with signature $\varepsilon$ is obtained by comparing the $\Q(\itPi)$-rational structures on both sides of the isomorphism.
We can arrange the periods $p({}^\sigma\!\itPi,\varepsilon)$ for $\sigma \in {\rm Aut}(\C)$ in a compatible way and obtain an element
\[
\left( p({}^\sigma\!\itPi,\varepsilon) \right)_{\sigma : \Q(\itPi)\rightarrow \C} \in (\Q(\itPi)\otimes_\Q\C)^\times/\Q(\itPi)^\times.
\]
Following is the main result of this paper, we prove a relation between the Betti--Whittaker periods of $\itPi$ and its contragredient $\itPi^\vee$.

\begin{thm}[Theorem \ref{T:main 1}]\label{T:Main 1}
Let $\itPi$ be a regular algebraic cuspidal automorphic representation of $\GL_n(\A)$.
Let $\varepsilon \in \{\pm1\}$ if $n$ is even, and $\varepsilon = \varepsilon(\itPi_\infty)$ if $n$ is odd.
Assume the regularity conditions in (\ref{E:regularity 1}) and (\ref{E:regularity 2}) below are satisfied.
Then we have 
\[
\sigma \left(\frac{p(\itPi,\varepsilon)}{G(\omega_\itPi)^{n-1}\cdot p(\itPi^\vee,\varepsilon)} \right) = \frac{p({}^\sigma\!\itPi,\varepsilon)}{G(\omega_{{}^\sigma\!\itPi})^{n-1}\cdot p({}^\sigma\!\itPi^\vee,\varepsilon)} ,\quad \sigma \in {\rm Aut}(\C).
\]
Here $G(\omega_\itPi)$ is the Gauss sum of the central character $\omega_\itPi$ of $\itPi$.
\end{thm}

\begin{rmk}
When $n=2$, Theorem \ref{T:Main 1} is compatible with the result of Raghuram and Shahidi \cite{RS2008} since $\itPi^\vee = \itPi \otimes \omega_\itPi^{-1}$.
\end{rmk}

\begin{rmk}
By Theorem \ref{T:Main 1} and the result of Balasubramanyam and Raghuram \cite{BR2017}, the analogous period relation also holds for the top degree Betti--Whittaker periods.
\end{rmk}

\begin{rmk}
In this paper, we work over $\Q$ for simplicity of exposition. 
It seems likely that the main theorem and its proof will work over a general number field, however we did not carry out the details in this paper.
One serious obstruction is the existence result Lemma \ref{L:auxiliary 1} which holds for totally real fields or CM-fields but is unknown for general number field.
\end{rmk}

As a consequence of the period relation, we prove the algebraicity of ratios of successive critical $L$-values for ${\rm GSpin}_{2n}^* \times \GL_{n'}$, where ${\rm GSpin}_{2n}^*$ is a quasi-split general spin group over $\Q$ of type $D_n$ or ${}^2D_n$. This is a generalization of the result of Bhagwat and Raghuram \cite{BR2020} for ${\rm O}_{n,n} \times \GL_1$ where $n$ is even and ${\rm O}_{n,n}$ is the even split orthogonal group of rank $n$ (cf.\,Remark \ref{R:BR} below).
We prove the algebraicity based on the period relation in Theorem \ref{T:Main 1} for orthogonal representations of $\GL_{2n}(\A)$ and the result of Harder and Raghuram \cite{HR2020} for $\GL_{2n} \times \GL_{n'}$, whereas in \cite{BR2020} the authors work intrinsically within even split orthogonal groups and study in detail the rank-one Eisenstein cohomology of ${\rm O}_{n+1,n+1}$.
A cuspidal automorphic representation $\itPi$ of $\GL_{2n}(\A)$ is called $\chi$-orthogonal for some Hecke character $\chi$ if the twisted symmetric square $L$-function $L(s,\itPi,{\rm Sym}^2\otimes\chi^{-1})$ has a pole at $s=1$.
In this case, $\itPi$ descends to a globally generic cuspidal automorphic representation of ${\rm GSpin}_{2n}^*(\A)$ and the general spin group is determined by the quadratic character $\chi^n\omega_\itPi^{-1}$ (cf.\,\cite{AS2006b}, \cite{AS2014}, and \cite{HS2016}).
Following is the second main result of this paper:

\begin{thm}[Theorem \ref{T:main 2}]\label{T:Main 2}
Let $\itPi$ and $\itSigma$ be regular algebraic cuspidal automorphic representations of $\GL_{2n}(\A)$ and $\GL_{n'}(\A)$ respectively. 
Let 
\[
L(s, \itPi \times \itSigma)
\]
be the associated Rankin--Selberg $L$-function.
Assume $\itPi$ is $\chi$-orthogonal and satisfies the regularity condition in (\ref{E:regularity 2}) below. Let $m_0,m_0+1 \in \Z+\tfrac{n'}{2}$ be critical points for $L(s,\itPi \times \itSigma)$ such that $L(m_0+1,\itPi \times \itSigma) \neq 0$. Then we have
\begin{align*}
&\sigma \left( \frac{L(m_0,\itPi \times \itSigma)}{(\sqrt{-1})^{nn'}\cdot G(\chi^{n}\omega_\itPi^{-1})^{n'}\cdot L(m_0+1,\itPi \times \itSigma)} \right)\\
& =  \frac{L(m_0,{}^\sigma\!\itPi\times{}^\sigma\!\itSigma)}{(\sqrt{-1})^{nn'}\cdot G({}^\sigma\!\chi^{n}\omega_{{}^\sigma\!\itPi}^{-1})^{n'}\cdot L(m_0+1,{}^\sigma\!\itPi\times{}^\sigma\!\itSigma)} ,\quad \sigma \in {\rm Aut}(\C).
\end{align*}
\end{thm}

\begin{rmk}\label{R:BR}
By the result of Bhagwat and Raghuram \cite{BR2020} for ${\rm O}_{n,n} \times \GL_1$, the assertion in Theorem \ref{T:Main 2} also holds under the following assumptions:
\begin{itemize}
\item If $L(s,\itPi)$ admits precisely $3$ critical points, then $n'=1$.
\item $n$ is even.
\item $\itPi$ is $|\mbox{ }|_\A$-orthogonal and $\omega_\itPi = |\mbox{ }|_\A^n$.
\end{itemize}
The last assumption is equivalent to saying that $\itPi\otimes|\mbox{ }|^{-1/2}_\A$ descends to a globally generic cuspidal automorphic representation of $\SO_{n,n}(\A)$.
\end{rmk}

Let $\F$ be an \'{e}tale real quadratic algebra over $\Q$. We denote by $\infty_1$ and $\infty_2$ the non-zero algebra homomorphisms from $\F$ into $\R$. Consider the Langlands functoriality of the Asai transfer from $\GL_{2}(\A_\F)$ to $\GL_{4}(\A)$. The Asai transfer of a regular algebraic cuspidal automorphic representation of $\GL_2(\A_\F)$ is either cuspidal or an isobaric sum of orthogonal cuspidal automorphic representations of $\GL_2(\A)$ with respect to a same Hecke character (cf.\,Corollaries \ref{C:tensor} and \ref{C:Asai}).
By Theorem \ref{T:Main 2}, we obtain the following result on the twisted Asai $L$-functions:

\begin{corollary}[Corollary \ref{C:main 2}]
Let $\pi$ and $\Sigma$ be regular algebraic cuspidal automorphic representations of $\GL_2(\A_\F)$ and $\GL_{n'}(\A)$ respectively. 
Let ${\rm As}(\pi)$ be the Asai transfer of $\pi$ to $\GL_4(\A)$.
Assume the following regularity condition is satisfied:
\[
\min\{\kappa_1,\kappa_2\}\geq \begin{cases}
3 & \mbox{ if $\kappa_1+\kappa_2$ is even},\\
4 & \mbox{ if $\kappa_1+\kappa_2$ is odd},
\end{cases}
\]
where $\kappa_i \geq 2$ is the minimal ${\rm SO}_2(\R)$-weight of $\pi_{\infty_i}$ for $i=1,2$.
Let $m_0,m_0+1 \in \Z+\tfrac{n'-1}{2}$ be critical points for $L(s,{\rm As}(\pi) \times \itSigma)$ such that $L(m_0+1, {\rm As}(\pi)\times \itSigma) \neq 0$. Then we have
\[
\sigma \left( \frac{L(m_0,{\rm As}(\pi) \times \itSigma)}{G(\omega_{\F/\Q})^{n'}\cdot L(m_0+1,{\rm As}(\pi) \times \itSigma)} \right)=  \frac{L(m_0, {\rm As}({}^\sigma\!\pi)\times {}^\sigma\!\itSigma)}{G(\omega_{\F/\Q})^{n'}\cdot L(m_0+1,{\rm As}({}^\sigma\!\pi)\times {}^\sigma\!\itSigma)} ,\quad \sigma \in {\rm Aut}(\C).
\]
Here $\omega_{\F/\Q}$ is the quadratic Hecke character associated to $\F/\Q$ by class field theory.
\end{corollary}

\begin{rmk}
The result is compatible with Deligne's conjecture by the period relation \cite[Theorem 3.6.2]{DR2023} due to Deligne and Raghuram. Please refer to \S\,\ref{SS:DR} below.
\end{rmk}

This paper is organized as follows. In \S\,\ref{S:BW}, we recall certain archimedean periods in \S\,\ref{SS:archimedean periods} and the Betti--Whittaker periods in \S\,\ref{SS:BW}. We stress that the Betti--Whittaker periods are defined under the normalization (\ref{E:normalization}) for generators in relative Lie algebra cohomology. In \S\,\ref{S:duality}, the main result is Theorem \ref{T:ratio under duality} on the algebraicity of ratios of critical values of Rankin--Selberg $L$-functions under duality which is the new input of this paper. In \S\,\ref{S:main results}, we recall orthogonal and symplectic representations in \S\,\ref{SS:polarized} and prove our main results in \S\,\ref{SS:main results}. In \S\,\ref{S:Deligne}, we show that our period relation in Theorem \ref{T:Main 1} is compatible with Deligne's conjecture on critical $L$-values of motives \cite[Conjecture 2.8]{Deligne1979}. More precisely, assume the validity of Clozel's and Deligne's conjectures, we show that Theorem \ref{T:Main 1} follows from the motivic period relation Lemma \ref{L:motivic}.

\subsection{Notations}\label{SS:notation}

Let $\A$ be the ring of adeles of $\Q$. Let $\A_f$ be the finite part of $\A$, and $\widehat{\Z} = \prod_p \Z_p$ be its maximal compact subring. 
For each place $v$ of $\Q$, let $|\mbox{ }|_v$ be the absolute value on $\Q_v$ normalized so that $|p|_p=p^{-1}$ if $v=p$ is a prime number and $|\mbox{ }|_\infty = |\mbox{ }|$ is the ordinary absolute value on $\R$ if $v = \infty$. 
Let $|\mbox{ }|_\A = \prod_v |\mbox{ }|_v$ be the normalized absolute value on $\A$.

Let $\chi$ be an algebraic Hecke character of $\A^\times$. We denote by $G(\chi)$ the Gauss sum of $\chi$ defined by
\[
G(\chi) = \prod_p\varepsilon(0,\chi_p,\psi_p),
\]
where $\varepsilon(s,\chi_p,\psi_p)$ is the $\varepsilon$-factor of $\chi_p$ with respect to $\psi_p$ defined in \cite{Tate1979}.
For $\sigma \in {\rm Aut}(\C)$, let ${}^\sigma\!\chi$ the unique algebraic Hecke character of $\A^\times$ such that ${}^\sigma\!\chi(a) = \sigma (\chi(a))$ for $a \in \A_f^\times$.
It is easy to verify that 
\begin{align}\label{E:Galois Gauss sum}
\begin{split}
\sigma(G(\chi)) = {}^\sigma\!\chi(u_\sigma)G({}^\sigma\!\chi),
\end{split}
\end{align} 
where $u_\sigma \in \widehat{\Z}^\times$ is the unique element such that $\sigma(\psi(x)) = \psi(u_\sigma x)$ for $x \in \A_f$.

Let $\psi=\bigotimes_v\psi_v$ be the standard additive character of $\Q\backslash \A$ defined so that
\begin{align*}
\psi_p(x) & = e^{-2\pi \sqrt{-1}x} \mbox{ for }x \in \Z[p^{-1}],\\
\psi_\infty(x) & = e^{2\pi \sqrt{-1}x} \mbox{ for }x \in \R.
\end{align*}
For $n \geq 1$, let $N_n$ be the standard maximal unipotent subgroup of $\GL_n$ consisting of upper unipotent matrices.
Let $\psi_{N_n} : N_n(\Q)\backslash N_n(\A) \rightarrow \C$ be the additive character defined by
\[
\psi_{N_n}(u) = \psi(u_{12}+u_{23}+\cdots+u_{n-1,n}),\quad u=(u_{ij}) \in N_n(\A).
\]
For each place $v$ of $\Q$, let $\psi_{N_n,v}$ be the local component of $\psi_{N_n}$ at $v$. 
Let $\psi_{N_n,f} = \bigotimes_{p} \psi_{N_n,p}$.

\section{Betti--Whittaker periods for $\GL_n$}\label{S:BW}

\subsection{Cohomological representations}

Let $K_n^\circ$ be the closed subgroup of $\GL_n(\R)$ defined by
\[
K_n^\circ = \R_+ \cdot \SO_n(\R).
\]
Here we embedded $\R_+$ into the center of $\GL_n(\R)$. We denote by $\frak{g}_n$ and $\frak{k}_n$ the Lie algebras of $\GL_n(\R)$ and $K_n$, respectively.
In this section, we recall some $(\frak{g}_n,K_n^\circ)$-cohomologcial representations of $\GL_n(\R)$ which are the archimedean local components of regular algebraic cuspidal automorphic representations of $\GL_n(\A)$.

Let $X^+(T_n)$ be the set of dominant integral weights for $\GL_n$ consisting of tuples of integers $\mu = (\mu_1,\cdots,\mu_n)$ with $\mu_1 \geq \cdots \geq \mu_n$. We say $\mu \in X^+(T_n)$ is pure if $\mu_{i}+\mu_{n+1-i} = \mu_{j}+\mu_{n+1-j}$ for all $1 \leq i,j \leq n$. Let $X_0^+(T_n)$ be the subset of $X^+(T_n)$ consisting of pure weights.
For $\mu \in X_0^+(T_n)$, let $(\underline{\kappa};\,{\sf w}) \in \Z^{\lfloor \tfrac{n}{2}\rfloor} \times \Z$ be the tuple of integers determined by ${\sf w} = -\mu_1-\mu_n$ and 
\begin{align*}
\mu = -\rho_n+
\begin{cases}
(\tfrac{\kappa_1-1-{\sf w}}{2},\cdots,\tfrac{\kappa_r-1-{\sf w}}{2},\tfrac{1-\kappa_r-{\sf w}}{2},\cdots,\tfrac{1-\kappa_1-{\sf w}}{2}) & \mbox{ if $n=2r$},\\
(\tfrac{\kappa_1-1-{\sf w}}{2},\cdots,\tfrac{\kappa_r-1-{\sf w}}{2},-\tfrac{{\sf w}}{2},\tfrac{1-\kappa_r-{\sf w}}{2},\cdots,\tfrac{1-\kappa_1-{\sf w}}{2}) & \mbox{ if $n=2r+1$}.
\end{cases}
\end{align*}
Here $\rho_n = (\tfrac{n-1}{2},\tfrac{n-3}{2},\cdots,\tfrac{1-n}{2})$ is half the sum of positive roots.
It is easy to verify that 
\begin{align}\label{E:infinity type}
\kappa_1>\cdots>\kappa_r \geq 2,\quad \begin{cases}
\kappa_i \equiv {\sf w}\,({\rm mod}\,2) & \mbox{ if $n=2r$},\\
\kappa_i \equiv {\sf w}+1 \equiv 1\,({\rm mod}\,2) & \mbox{ if $n=2r+1$}.
\end{cases}
\end{align}
The association $\mu \leftrightarrow(\underline{\kappa};\,{\sf w})$ is a one-to-one correspondence between $X_0^+(T_n)$ and the subset of tuples in $\Z^{\lfloor \tfrac{n}{2}\rfloor} \times \Z$ satisfying (\ref{E:infinity type}).
For $\mu \in X_0^+(T_n)$, let $\pi_\mu$ be the irreducible admissible $(\frak{g}_n,{\rm O}_n(\R))$-module realized as the space of ${\rm O}_n(\R)$-finite vectors of the following induced representation of $\GL_n(\R)$: 
\begin{align*}
\begin{cases}
{\rm Ind}_{P_{(2,\cdots,2)}(\R)}^{\GL_n(\R)}\left(D_{\kappa_{1}} {\otimes}\cdots {\otimes} D_{\kappa_{r}}\right)\otimes |\mbox{ }|^{{\sf w}/2} & \mbox{ if $n=2r$},\\
{\rm Ind}_{P_{(2,\cdots,2,1)}(\R)}^{\GL_n(\R)}\left(D_{\kappa_{1}} {\otimes}\cdots {\otimes} D_{\kappa_{r}}\otimes 1\right)\otimes |\mbox{ }|^{{\sf w}/2} & \mbox{ if $n=2r+1$}.
\end{cases}
\end{align*}
Here $P_{(n_1,\cdots,n_k)}$ is the standard upper parabolic subgroup of $\GL_n$ of type $(n_1,\cdots,n_k)$, and $D_\kappa$ is the discrete series representation of $\GL_2(\R)$ with minimal weight $\kappa \geq 2$. Let $\Omega_0(n)$ be the set of irreducible admissible $(\frak{g}_n,{\rm O}_n(\R))$-modules defined by
\[
\Omega_0(n) = \bigcup_{\mu \in X_0^+(T_n)}\Omega_\mu,
\]
where 
\[
\Omega_\mu = \begin{cases}
\left\{\pi_\mu\right\} & \mbox{ if $n$ is even},\\
\left\{\pi_\mu,\pi_\mu\otimes{\rm sgn}\right\} & \mbox{ if $n$ is odd}.
\end{cases}
\] 

Let $\itPi \in \Omega_\mu$. We call $(\underline{\kappa};\,{\sf w})$ the infinity type of $\itPi$.
When $n=2r+1$, let $\varepsilon(\itPi) \in \{\pm1\}$ be the signature of $\itPi$ defined by 
\begin{align}\label{E:signature}
\varepsilon(\itPi) = \begin{cases}
(-1)^{r+{\sf w}/2} & \mbox{ if $\itPi = \pi_\mu$},\\
(-1)^{r+1+{\sf w}/2} & \mbox{ if $\itPi = \pi_\mu\otimes{\rm sgn}$}.
\end{cases}
\end{align}
By \cite[Theorem 3.3, III]{BW2000}, $\itPi$ is $(\frak{g}_n,K_n^\circ)$-cohomological with coefficients in $M_{\mu,\C}$, that is, we have (cf.\,\cite[Lemme 3.14]{Clozel1990})
\[
H^\bullet(\frak{g}_n,K_n^\circ;\itPi\otimes M_{\mu,\C}) \neq 0.
\]
The group $\pi_0(\GL_n(\R)) = \GL_n(\R) / \GL_n(\R)^\circ \cong \Z/2\Z$ naturally acts on these cohomology groups and we denote by $H^\bullet(\frak{g}_n,K_n^\circ;\itPi\otimes M_{\mu,\C})[\varepsilon]$ the $\varepsilon$-isotypic component under the action for $\varepsilon \in \{\pm1\}$.
In particular, when
\[
\bullet = b_n = \lfloor \tfrac{n^2}{4} \rfloor
\]
is the bottom degree, we have
\[
{\rm dim}_\C\,H^{b_n}(\frak{g}_n,K_n^\circ;\itPi\otimes M_{\mu,\C})[\varepsilon] = \begin{cases}
1 & \mbox{ if $n$ is even or $\varepsilon = \varepsilon(\itPi)$ if $n$ is odd},\\
0 & \mbox{ otherwise}.
\end{cases}
\]
From now on, we assume $\varepsilon = \varepsilon(\itPi)$ if $n$ is odd and fix a generator
\begin{align}\label{E:generator}
[\itPi]^\varepsilon \in H^{b_n}(\frak{g}_n,K_n^\circ;\mathcal{W}(\itPi)\otimes M_{\mu,\C})[\varepsilon],
\end{align}
where $\mathcal{W}(\itPi)$ is the Whittaker model of $\itPi$ with respect to $\psi_{N_n,\infty}$.
We then have a set of generators
\[
\left\{ [\itPi\otimes\chi]^{\varepsilon\cdot\varepsilon(\chi)}  \,\left\vert \, \itPi \in \Omega_0(n),\,\chi \in \Omega_0(1) \right\}\right.
\]
We normalize this set as follows: Let $\itPi \in \Omega_\mu\subset\Omega_0(n)$ and $\chi \in \Omega_0(1)$. We have $\chi = {\rm sgn}^\delta|\mbox{ }|^{{\sf u}}$ for some $\delta \in \{0,1\}$ and ${\sf u} \in \Z$. Note that $M_{\mu-{\sf u}} = M_\mu\otimes \det^{-{\sf u}}$ and we identify the representation spaces of $M_\mu$ and $M_{\mu-{\sf u}}$.
The $\C$-linear isomorphism 
\[
\mathcal{W}(\itPi) \longrightarrow \mathcal{W}(\itPi\otimes\chi),\quad W \longmapsto W\otimes\chi
\]
and the identity map $M_{\mu} \rightarrow M_{\mu-{\sf u}}$ then induces isomorphism of complexes for relative Lie algebra cohomology 
\[
\left( \text{\large$ \wedge $}^\bullet(\frak{g}_{n,\C}/\frak{k}_{n,\C})^*\otimes \mathcal{W}(\itPi) \otimes M_{\mu,\C} \right)^{K_n^\circ} \longrightarrow \left( \text{\large$ \wedge $}^\bullet(\frak{g}_{n,\C}/\frak{k}_{n,\C})^* \otimes \mathcal{W}(\itPi\otimes\chi) \otimes M_{\mu-{\sf u},\C}\right)^{K_n^\circ},
\] 
which in turn defines $\C$-linear isomorphisms
\[
A_\chi^\bullet:H^\bullet(\frak{g}_n,K_n^\circ;\mathcal{W}(\itPi)\otimes M_{\mu,\C}) \longrightarrow H^\bullet(\frak{g}_n,K_n^\circ;\mathcal{W}(\itPi\otimes\chi)\otimes M_{\mu-{\sf u},\C}).
\]
In particular, when $\bullet = b_n$, $A_\chi^{b_n}$ will send $[\itPi]^\varepsilon$ to a scalar multiple of $[\itPi\otimes\chi]^{\varepsilon\cdot\varepsilon(\chi)}$.
We normalize these generators so that 
\begin{align}\label{E:normalization}
A_\chi^{b_n}([\itPi]^\varepsilon) = [\itPi\otimes\chi]^{\varepsilon\cdot\varepsilon(\chi)}.
\end{align}

\subsection{Archimedean periods for $\GL_{n} \times \GL_{n-1}$}\label{SS:archimedean periods}

In this section, we recall certain archimedean periods for $\GL_{n}(\R) \times \GL_{n-1}(\R)$ which appeared in the algebraicity result of Raghuram \cite{Raghuram2009} for $\GL_{n}(\A) \times \GL_{n-1}(\A)$.
Let $\itPi \in \Omega_{\mu}\subset\Omega_0(n)$ and $\itSigma \in \Omega_{\lambda}\subset\Omega_0(n-1)$ with infinity types  $(\underline{\kappa};\,{\sf w})$ and $(\underline{\ell};\,{\sf u})$ respectively. We say $(\itPi,\itSigma)$ is balanced if
\begin{align}\label{E:balanced}
\begin{cases}
\kappa_1> \ell_1 > \kappa_2 > \ell_2>\cdots>\kappa_{r-1}>\ell_{r-1}>\kappa_r & \mbox{ if $n=2r$},\\
\kappa_1> \ell_1 > \kappa_2 > \ell_2>\cdots>\kappa_{r-1}>\ell_{r-1}>\kappa_r>\ell_r& \mbox{ if $n=2r+1$}.
\end{cases}
\end{align}
Assume $(\itPi,\itSigma)$ is balanced. By the branching law from $\GL_{n1}$ to $\GL_{n-1}$, the Rankin--Selberg $L$-factors (cf.\,\S\,\ref{SS:RS}) $L(s,\itPi\times\itSigma)$ and $L(1-s,\itPi^\vee \times\itSigma^\vee)$ are holomorphic at a half-integer $m+\tfrac{1}{2} \in \Z+\tfrac{1}{2}$ if and only if 
\[
{\rm Hom}_{\GL_{n-1}(\Q)}(M_{\mu}\otimes M_{\lambda},\mbox{$\det^m$}) \neq 0.
\]
In this case, the above space is one-dimensional and we fix a non-zero functional $f_{\mu,\lambda,m}$. 
Let 
\[
Z_m(\cdot,\cdot): H^{b_{n}}(\frak{g}_{n},K_{n}^\circ;\mathcal{W}(\itPi) \otimes M_{\mu,\C}) \times H^{b_{n-1}}(\frak{g}_{n-1},K_{n-1}^\circ;\mathcal{W}(\itSigma) \otimes M_{\lambda,\C}) \longrightarrow \C
\]
be a pairing defined as follows: 
For $W \in \mathcal{W}(\itPi)$ and $W' \in \mathcal{W}(\itSigma)$, let $Z(s,W,W')$ be the local zeta integral defined by
\[
Z(s,W,W') = \int_{N_{n-1}(\R)\backslash \GL_{n-1}(\R)}W\left( \bp g & 0 \\ 0 & 1 \ep \right)W'(\delta_{n-1}g)|\det g|^{s-1/2}\,dg,\\
\]
where $\delta_{n-1} = {\rm diag}(-1,1,\cdots,(-1)^{n-1})$. Here $dg$ is the quotient measure of the Haar measure 
\[
\prod_{i=1}^{n-1}\Gamma_\R(i)\cdot \frac{\prod_{1 \leq i,j \leq {n-1}}dg_{ij}}{|\det g|^{n-1}}
\]
on $\GL_{n-1}(\R)$ by the Haar measure on $N_{n-1}(\R)$ given by the product of Lebesgue measures.
The integral converges absolutely for ${\rm Re}(s)$ sufficiently large and admits meromorphic continuation to $s \in \C$.
Moreover, it is holomorphic at $s=m+\tfrac{1}{2}$ (cf.\,\cite[Theorem 2.1-(ii)]{Jacquet2009}). 
Let $e_{ij} \in \frak{g}_{n-1}$ be the $(n-1) \times (n-1)$ matrix with $1$ in the $(i,j)$-entry and zeros otherwise and $\{e_{ij}^*\,\vert\,1 \leq i,j\leq n-1\}$ be the corresponding dual basis of $\frak{g}_{n-1}^*$ ordered lexicographically.
Let 
\[
(\cdot,\cdot) : \text{\large$ \wedge $}^{b_{n}}(\frak{g}_{n,\C}/\frak{k}_{n,\C})^* \times \text{\large$ \wedge $}^{b_{n-1}}(\frak{g}_{n-1,\C}/\frak{k}_{n-1,\C})^* \longrightarrow \C
\]
be the pairing defined so that
\[
\iota(X^*)\wedge {\rm pr}(Y^*) = (X^*,Y^*)\cdot \wedge_{1 \leq i\leq j \leq n-1}e_{ij}^*,
\]
where $\iota : \text{\large$ \wedge $}^{b_{n}}(\frak{g}_{n,\C}/\frak{k}_{n,\C})^* \rightarrow \text{\large$ \wedge $}^{b_{n}}(\frak{g}_{n-1,\C}/\frak{so}_{n-1,\C})^*$ is induced from the inclusion $\GL_{n-1}(\R) \subset \GL_{n}(\R)$ sending $g$ to ${\rm diag}(g,1)$, and ${\rm pr}: \text{\large$ \wedge $}^{b_{n-1}}(\frak{g}_{n-1,\C}/\frak{k}_{n-1,\C})^* \rightarrow \text{\large$ \wedge $}^{b_{n-1}}(\frak{g}_{n-1,\C}/\frak{so}_{n-1,\C})^*$ is the natural projection.
We then define $Z_m(\cdot,\cdot)$ to be the restriction of 
\[
(\cdot,\cdot)\otimes Z(s,\cdot,\cdot)\vert_{s=m+{1}/{2}} \otimes f_{\mu,\lambda,m}.
\]
Note that the pairing $Z_m(\cdot,\cdot)$ is well-defined up to scalar multiple of $\Q^\times$ as it depends on the choice of $f_{\mu,\lambda,m}$.
We recall the following important non-vanishing result of Sun \cite{Sun2017}.

\begin{thm}[Sun]\label{T:Sun}
Let $\itPi \in \Omega_0(n)$ and $\itSigma \in \Omega_0(n-1)$. Assume $(\itPi,\itSigma)$ is balanced. 
For $m\in\Z$ such that $L(s,\itPi\times\itSigma)$ and $L(1-s,\itPi^\vee\times\itSigma^\vee)$ are holomorphic at $s=m+\tfrac{1}{2}$, we have
\[
Z_m([\itPi]^{\varepsilon_m},[\itSigma]^{\varepsilon_m'}) \neq 0.
\]
Here $\varepsilon_m = \varepsilon(\itPi)$ if $n$ is odd, $\varepsilon_m' = \varepsilon(\itSigma)$ if $n$ is even, and $\varepsilon_m\varepsilon_m' = (-1)^{m+n}$.
\end{thm}
\begin{definition}\label{D:archimedean periods}
Let notation and assumption be as in Theorem \ref{T:Sun}.
We define the archimedean period $p(m,\itPi \times \itSigma) \in \C^\times/\Q^\times$ by
\[
p(m,\itPi \times \itSigma) = \frac{L(m+\tfrac{1}{2},\itPi \times \itSigma)}{Z_m([\itPi]^{\varepsilon_m},[\itSigma]^{\varepsilon_m'}) }.
\]
\end{definition}

The following lemma is on the relation between the archimedean periods upon twisting by powers of the absolute value.

\begin{lemma}\label{E:archi. period comparison}
Let $\itPi \in \Omega_0(n)$ and $\itSigma \in \Omega_0(n-1)$. Assume $(\itPi,\itSigma)$ is balanced. For ${\sf w}_1,{\sf w}_2 \in \Z$,  we have
\[
\frac{p(m,(\itPi \otimes |\mbox{ }|^{{\sf w}_1}) \times (\itSigma\otimes|\mbox{ }|^{{\sf w}_2}))}{p(m+{\sf w}_1+{\sf w}_2,\itPi \times \itSigma)} \in \Q^\times
\]
for all $m \in \Z$ such that $L(s,\itPi \times \itSigma)$ and $L(1-s,\itPi^\vee \times \itSigma^\vee)$ are holomorphic at $s = m+{\sf w}_1+{\sf w}_2+\tfrac{1}{2}$.
\end{lemma}

\begin{proof}
Assume $\itPi \in \Omega_\mu$ and $\itSigma \in \Omega_\lambda$.
Put $m' = m+{\sf w}_1+{\sf w}_2$.
For $W \in \mathcal{W}(\itPi)$ and $W' \in \mathcal{W}(\itSigma)$, it is clear that
\[
Z(s+{\sf w}_1+{\sf w}_2,W,W') = Z(s,W\otimes |\mbox{ }|^{{\sf w}_1},W'\otimes |\mbox{ }|^{{\sf w}_2}).
\]
The identity maps $M_{\mu} \rightarrow M_{\mu-{\sf w}_1}$ and $M_{\lambda} \rightarrow M_{\lambda-{\sf w}_2}$ induce an isomorphism 
\[
{\rm Hom}_{\GL_{n-1}(\Q)}(M_{\mu}\otimes M_{\lambda},\mbox{$\det^{m'}$}) \longrightarrow {\rm Hom}_{\GL_{n-1}(\Q)}(M_{\mu-{\sf w}_1}\otimes M_{\lambda-{\sf w}_2},\mbox{$\det^m$}),
\] 
which sends $f_{\mu,\lambda,m'}$ to $C\cdot f_{\mu-{\sf w}_1,\lambda-{\sf w}_2,m}$ for some $C \in \Q^\times$.
Therefore, we have
\begin{align*}
Z_{m'}([\itPi]^{\varepsilon_{m'}},[\itSigma]^{\varepsilon'_{m'}}) &= C\cdot Z_{m}(A_{|\mbox{ }|^{{\sf w}_1}}^{b_{n}}([\itPi]^{\varepsilon_{m'}}),A_{|\mbox{ }|^{{\sf w}_2}}^{b_{n-1}}([\itSigma]^{\varepsilon_{m'}'}))\\
& = C\cdot Z_{m}([\itPi\otimes|\mbox{ }|^{{\sf w}_1}]^{\varepsilon_{m}},[\itSigma\otimes|\mbox{ }|^{{\sf w}_2}]^{\varepsilon'_{m}}).
\end{align*}
Here the second equality follows from our normalization (\ref{E:normalization}). 
This completes the proof.
\end{proof}

\subsection{Betti--Whittaker periods}\label{SS:BW}

Consider the topological space
\[
\mathcal{S}_n = \GL_n(\Q)\backslash\GL_n(\A)/K_n^\circ.
\]
For $\mu \in X^+(T_n)$, let $\mathcal{M}_\mu$ be the sheaf of $\Q$-vector spaces on $\mathcal{S}_n$ associated to $M_\mu$ (cf.\,\cite[\S\,2.2.8]{HR2020}). 
The sheaf cohomology group of $\mathcal{M}_\mu$ on $\mathcal{S}_n$ is denoted by
\[
H^\bullet(\mathcal{S}_n,\mathcal{M}_\mu).
\]
The group is naturally equipped with action of $\pi_0(\GL_n(\R))\times\GL_n(\A_f)$.
At the transcendental level, we have a canonical isomorphism
\[
H^\bullet(\mathcal{S}_n,\mathcal{M}_{\mu,\C}) \cong H^\bullet(\frak{g}_{n},K_n^\circ;C^\infty(\GL_n(\Q)\backslash\GL_n(\A),\xi_\mu)\otimes M_{\mu,\C}),
\]
where $\xi_\mu$ is the character of $\R_+$ such that $\R_+$ acts on $M_{\mu,\C}$ by $\xi_\mu^{-1}$, and
$C^\infty(\GL_n(\Q)\backslash\GL_n(\A),\xi_\mu)$ is the space of smooth functions $\varphi$ on $\GL_n(\Q)\backslash\GL_n(\A)$ such that $\varphi(ga) = \xi_\mu(a)\varphi(g)$ for all $a \in \R_+$.
The cuspidal cohomology group of $\GL_n$ with coefficients in $M_{\mu,\C}$ is defined by
\[
H^\bullet_{\rm cusp}(\mathcal{S}_n,\mathcal{M}_{\mu,\C}) = H^\bullet(\frak{g}_{n},K_n^\circ;\mathcal{A}_0(\GL_n(\A),\xi_\mu)\otimes M_{\mu,\C}),
\]
where $\mathcal{A}_0(\GL_n(\A),\xi_\mu)$ is the subspace of $C^\infty(\GL_n(\Q)\backslash\GL_n(\A),\xi_\mu)$ consisting of smooth cusp forms on $\GL_n(\A)$.
The natural inclusion then induces a $(\pi_0(\GL_n(\R))\times\GL_n(\A_f))$-equivariant homomorphism 
\[
H^\bullet_{\rm cusp}(\mathcal{S}_n,\mathcal{M}_{\mu,\C}) \longrightarrow H^\bullet(\mathcal{S}_n,\mathcal{M}_{\mu,\C})
\]
which is injective by the result of Borel \cite[Corollary 5.5]{Borel1981}.
For $\sigma \in {\rm Aut}(\C)$, the $\sigma$-linear isomorphism 
\[
M_{\mu,\C} = M_\mu\otimes\C \longrightarrow M_{\mu,\C},\quad {\bf v}\otimes z \longmapsto {\bf v}\otimes\sigma(z)
\]
induces a $(\pi_0(\GL_n(\R))\times\GL_n(\A_f))$-equivariant $\sigma$-linear isomorphism 
\[
\sigma^\bullet : H^\bullet(\mathcal{S}_n,\mathcal{M}_{\mu,\C}) \longrightarrow H^\bullet(\mathcal{S}_n,\mathcal{M}_{\mu,\C}).
\]
By the results of Clozel \cite[Th\'{e}or\`{e}me 3.19 and Lemme 4.9]{Clozel1990}, the cuspidal cohomology is non-zero only if $\mu \in X_0^+(T_n)$ and we have
\begin{align}\label{E:Galois invariant}
\sigma^\bullet(H^\bullet_{\rm cusp}(\mathcal{S}_n,\mathcal{M}_{\mu,\C}) ) =H^\bullet_{\rm cusp}(\mathcal{S}_n,\mathcal{M}_{\mu,\C}) ,\quad \sigma \in {\rm Aut}(\C). 
\end{align}
Assume $\mu \in X_0^+(T_n)$. Let $\itPi$ be a cuspial automorphic representation of $\GL_n(\A)$ such that $\itPi_f = \bigotimes_p\itPi_p$ appears in the cuspidal cohomology $H^\bullet_{\rm cusp}(\mathcal{S}_n,\mathcal{M}_{\mu,\C})$, that is, $\itPi$ is regular algebraic and $\itPi_\infty \in \Omega_\mu$.
For $\sigma \in {\rm Aut}(\C)$, the $\sigma$-conjugate ${}^\sigma\!\itPi_f$ of $\itPi_f$ also appears in the cuspidal cohomology by (\ref{E:Galois invariant}). Hence 
\[
{}^\sigma\!\itPi=\itPi_\infty\otimes{}^\sigma\!\itPi_f
\]
is cuspidal automorphic.
Denote by $\Q(\itPi)$ the rationality field of $\itPi$, which is defined to be the fixed field of $\{\sigma \in {\rm Aut}(\C)\,\vert\,{}^\sigma\!\itPi = \itPi\}$.
It is a number field by (\ref{E:Galois invariant}) and the admissibility of $H^\bullet(\mathcal{S}_n,\mathcal{M}_{\mu,\C})$.
By abuse of notation, we write $\itPi$ for the representation space of $\itPi$ realized in $\mathcal{A}_0(\GL_n(\A))$.
For $\varphi \in \itPi$, let $W(\varphi)$ be the Whittaker function of $\varphi$ with respect to $\psi_{N_n}$ defined by
\[
W(g,\varphi) = \int_{N_n(\Q)\backslash N_n(\A)}\varphi(ug)\overline{\psi_{N_n}(u)}\,du^{\rm Tam},\quad g \in \GL_n(\A).
\]
Here $du^{\rm Tam}$ is the Tamagawa measure on $N_n(\A)$.
Let $\mathcal{W}(\itPi)$ be the space of Whittaker functions of $\itPi$. For each place $v$ of $\Q$, let $\mathcal{W}(\itPi_v)$ be the space of Whittaker functions of $\itPi_v$ with respect to $\psi_{N_n,v}$. 
Note that when $n=1$ and $\itPi=\chi$ is an algebraic Hecke character, we understand $\mathcal{W}(\chi) = \C\cdot\chi$ and $\mathcal{W}(\chi_v)=\C\cdot\chi_v$.
It is well-known that we have an isomorphism
\[
\bigotimes_v{'}\, \mathcal{W}(\itPi_v) \longrightarrow \mathcal{W}(\itPi),\quad \bigotimes_vW_v \longmapsto \prod_v W_v,
\]
where the restricted tensor product is defined with respect to the $\GL_n(\Z_p)$-invariant Whittaker function $W_{\itPi_p}^\circ \in \mathcal{W}(\itPi_p)$ normalized so that $W_{\itPi_p}^\circ({\bf 1}_n)=1$ for all primes $p$ at which $\itPi$ is unramified.
We denote by $\Upsilon_\itPi$ the inverse of the isomorphism $\itPi \rightarrow \mathcal{W}(\itPi)$ sending $\varphi$ to $W(\varphi)$. It then induces a $\pi_0(\GL_n(\R))$-equivariant injective homomorphism
\[
\Upsilon_\itPi^\bullet : H^\bullet(\frak{g}_{n},K_n^\circ;\mathcal{W}(\itPi)\otimes M_{\mu,\C}) \longrightarrow H^\bullet_{\rm cusp}(\mathcal{S}_n,\mathcal{M}_{\mu,\C}).
\]
Let $\varepsilon \in \{\pm1\}$ if $n$ is even, and $\varepsilon = \varepsilon(\itPi_\infty)$ if $n$ is odd.
Let 
\[
\itPhi_\itPi^\varepsilon :  \mathcal{W}(\itPi_f) = \bigotimes_p{'}\,\mathcal{W}(\itPi_p) \longrightarrow H^{b_n}_{\rm cusp}(\mathcal{S}_n,\mathcal{M}_{\mu,\C})
\]
be the injective homomorphism defined by
\[
\itPhi_\itPi^\varepsilon = \Upsilon_\itPi^{b_n}\circ ([\itPi_\infty]^\varepsilon\otimes\,\cdot\,),
\]
where $[\itPi_\infty]^\varepsilon$ is the generator fixed in 
(\ref{E:generator}).
For $\sigma \in {\rm Aut}(\C)$, let $t_\sigma : \mathcal{W}(\itPi_f) \rightarrow \mathcal{W}({}^\sigma\!\itPi_f)$ be the $\GL_n(\A_f)$-equivariant $\sigma$-linear isomorphism defined by
\[
t_\sigma W(g) = \sigma \left( W({\rm diag}(u_\sigma^{1-n},u_\sigma^{2-n},\cdots,1)g) \right),\quad g \in \GL_n(\A_f),
\]
where $u_\sigma \in \widehat{\Z}^\times$ is the finite idele such that $\sigma(\psi(x)) = \psi(u_\sigma x)$ for $x \in \A_f$.
By comparing the $\Q(\itPi)$-rational structures on $\itPi_f$ given by $\mathcal{W}(\itPi_f)$ via the action of $t_\sigma$ and by the $(\varepsilon\times\itPi_f)$-isotypic component in the cuspidal cohomology of degree $b_n$, we have the following definition of the (bottom degree) Betti--Whittaker period introduced by Mahnkopf \cite{Mahnkopf2005} (see also \cite{Harder1983}, \cite{Hida1994} ($n=2$) and \cite{RS2008} for general number field).

\begin{definition}
Let $\itPi$ be a regular algebraic cuspidal automorphic representation of $\GL_n(\A)$.
Let $\varepsilon \in \{\pm1\}$ if $n$ is even, and $\varepsilon = \varepsilon(\itPi_\infty)$ if $n$ is odd.
Under the canonical isomorphism $\Q(\itPi)\otimes_\Q\C\cong \prod_{\sigma : \Q(\itPi)\rightarrow \C}\C$, there exists a unique element 
\[
(p({}^\sigma\!\itPi,\varepsilon))_{\sigma : \Q(\itPi)\rightarrow \C} \in (\Q(\itPi)\otimes_\Q\C)^\times / \Q(\itPi)^\times
\]
such that
\[
\sigma^{b_n}\circ\left( \frac{\Phi_\itPi^\varepsilon}{p(\itPi,\varepsilon)} \right) = \left( \frac{\Phi_{{}^\sigma\!\itPi}^\varepsilon}{p({}^\sigma\!\itPi,\varepsilon)} \right)\circ t_\sigma,\quad \sigma \in {\rm Aut}(\C).
\]
We call $p(\itPi,\varepsilon) \in \C^\times/\Q(\itPi)^\times$ the Betti--Whittaker period of $\itPi$ with signature $\varepsilon$.
\end{definition}

\begin{rmk}
If we replace $[\itPi_\infty]^\varepsilon$ by $C\cdot [\itPi_\infty]^\varepsilon$ for some $C\in\C^\times$, then the periods are replaced by
\[
(1\otimes C)\cdot (p({}^\sigma\!\itPi,\varepsilon))_{\sigma : \Q(\itPi)\rightarrow \C}.
\]
\end{rmk}

\begin{rmk}\label{R:n=1}
When $n=1$, $\itPi=\chi$ is an algebraic Hecke character. In this case, $b_{0}=0$ and we are reduced to consider the global sections of $\mathcal{M}_\mu$ on $\mathcal{S}_1$ at the transcendental level. Then it is clear that (cf.\,\cite[Lemma 1.25]{GL2021})
\[
(p({}^\sigma\!\chi,\varepsilon))_{\sigma : \Q(\chi)\rightarrow \C} \in (1\otimes C)\cdot\Q(\chi)^\times,
\]
where $C \in \C^\times$ is determined by
\[
[\chi_\infty]^\varepsilon = C\cdot (1 \otimes \chi_\infty \otimes 1) \in H^0(\frak{g}_{1,\C},\frak{k}_{1,\C};\mathcal{W}(\chi_\infty)\otimes M_{\mu,\C})[\varepsilon].
\]
\end{rmk}




\section{Ratios of Rankin--Selberg $L$-functions under duality}\label{S:duality}

The main result of this section is Theorem \ref{T:ratio under duality}.
We prove the algebraicity of ratios of critical values of Rankin--Selberg $L$-functions under duality which is a crucial ingredient in the proof of Theorem \ref{T:main 1}.
The result itself can be regarded as an automorphic analogue of the period relation in \cite[Proposition 5.1]{Deligne1979} for tensor product motives.

\subsection{Rankin--Selberg $L$-functions}\label{SS:RS}

Let $\itPi$ and $\itSigma$ be cuspidal automorphic representations of $\GL_{n}(\A)$ and $\GL_{n'}(\A)$ respectively.
For each place $v$ of $\Q$, let $\phi_{\itPi_v}$ and $\phi_{\itSigma_v}$ be the Langlands parameters of $\itPi_v$ and $\itSigma_v$ respectively. Associated to the tensor representation $\phi_{\itPi_v}\otimes\phi_{\itSigma_v}$ of the Weil--Deligne group of $\Q_v$, we have the local $L$-factor and $\varepsilon$-factor (cf.\,\cite[\S\,31.3]{BH2011})
\[
L(s,\itPi_v\times\itSigma_v) = L(s,\phi_{\itPi_v}\otimes\phi_{\itSigma_v}),\quad \varepsilon(s,\itPi_v \times \itSigma_v,\psi_v) = \varepsilon(s,\phi_{\itPi_v}\otimes\phi_{\itSigma_v},\psi_v).
\]
The Rankin--Selberg $L$-function associated to $\itPi$ and $\itSigma$ is defined by the Euler product
\[
L(s,\itPi\times\Sigma) = \prod_v L(s,\itPi_v\times\itSigma_v).
\]
It converges absolutely when ${\rm Re}(s)$ is sufficiently large and admits meromorphic continuation to the whole complex plane.
Moreover, we have the global functional equation
\begin{align}\label{E:global fe}
L(s,\itPi \times \itSigma) = \varepsilon(s,\itPi \times \itSigma) \cdot L(1-s,\itPi^\vee \times \itSigma^\vee),
\end{align}
where $\varepsilon(s,\itPi \times \itSigma) = \prod_v\varepsilon(s,\itPi_v \times \itSigma_v,\psi_v)$.
Assume further that $\itPi$ and $\itSigma$ are regular algebraic.
In this case, we have $\itPi \in \Omega(n)$ and $\itSigma \in \Omega(n')$.
A critical point is a half integer $m_0 \in \Z+\tfrac{n+n'}{2}$ such that $L(s,\itPi_\infty\times\itSigma_\infty)$ and $L(1-s,\itPi_\infty^\vee\times\itSigma_\infty^\vee)$ are holomorphic at $s=m_0$.
For instance, let $(\underline{\kappa};\,{\sf w})$ and $(\underline{\ell};\,{\sf u})$ be the infinity types of $\itPi_\infty$ and $\itSigma_\infty$ respectively. 
If $n$ is even, then the set of critical points is given by
\begin{align}\label{E:critical range}
\left \{ m_0 \in \Z+\tfrac{n'}{2} \,\left\vert\, \tfrac{2-{\sf w}-{\sf u}-d(\underline{\kappa},\underline{\ell})}{2} \leq m_0 \leq \tfrac{-{\sf w}-{\sf u}+d(\underline{\kappa},\underline{\ell})}{2} \right\}\right.,
\end{align}
where 
\[
d(\underline{\kappa},\underline{\ell}) = \begin{cases}
\min\{|\kappa_i-\ell_j|\} & \mbox{ if $n'$ is even},\\
\min\{|\kappa_i-\ell_j|,|\kappa_i-1|\} & \mbox{ if $n'$ is odd}.
\end{cases}
\]
In particular, $L(s,\itPi\times\itSigma)$ must be entire if it admits critical points. 
The central point $s = \tfrac{1-{\sf w}-{\sf u}}{2}$ is critical if and only if $d(\underline{\kappa},\underline{\ell})\geq 1$ and ${\sf w}+{\sf u}\equiv n+n'+1\,({\rm mod}\,2)$.
By the results of Jacquet--Shalika \cite[Theorem 5.3]{JS1981} and Shahidi \cite[Theorem 5.2]{Shahidi1981}, $L(s,\itPi\times\itSigma)$ is non-zero at non-central critical points.

\subsection{Galois equivariance of local factors}

In this section, we prove the Galois equivariance properties of the local Rankin--Selberg factors.
We interpreted the local $L$-factor as the greatest common divisor of the local zeta integrals introduced and studied by Jacquet, Piatetski-Shapiro, and Shalika in \cite{JS1981}, \cite{JS1981b}, and \cite{JPSS1983}.
The local $\varepsilon$-factor is also interpreted as the ratio appeared in the local functional equation.
First we recall the local zeta integrals. 
Fix a prime number $p$ in this section.
Let $\itPi_p$ and $\itSigma_p$ be irreducible admissible generic representations of $\GL_{n}(\Q_p)$ and $\GL_{n'}(\Q_p)$ respectively. 
Let $\mathcal{W}(\itPi_p)$ and $\mathcal{W}(\itSigma_p)$ be the Whittaker models of $\itPi_p$ and $\itSigma_p$ with respect to $\psi_{N_{n},p}$ and $\psi_{N_{n'},p}$ respectively.
Let $S(\Q_p^n)$ be the space of locally constant functions with compact supports on $\Q_p^n$.
For $W \in \mathcal{W}(\itPi_p)$, $W' \in \mathcal{W}(\itSigma_p)$, and $\Phi \in S(\Q_p^n)$, we define the local zeta integrals as follows:
If $n>n'$, let
\begin{align*}
Z(s,W,W') &= \int_{N_{n'}(\Q_p)\backslash \GL_{n'}(\Q_p)}W\left( \bp g & 0 \\ 0 & {\bf 1}_{n-n'} \ep \right)W'(\delta_{n'}g)|\det g|_p^{s-(n-n')/2}\,dg,\\
Z^\vee(s,W,W') &= \int_{N_{n'}(\Q_p)\backslash \GL_{n'}(\Q_p)} \int_{M_{n-n'-1,n'}(\Q_p)}W\left( \bp g & 0 & 0 \\x & {\bf 1}_{n-n'-1} & 0 \\ 0 & 0 & 1\ep \right)W'(\delta_{n'}g)\\
&\quad\quad\quad\quad\quad\quad\quad\quad\quad\quad\quad\quad\quad\quad\quad\quad\quad\quad\quad\quad|\det g|_p^{s-(n-n')/2}\,dx\,dg.
\end{align*}
If $n=n'$, let
\begin{align*}
Z(s,W,W',\Phi) = \int_{N_{n}(\Q_p)\backslash \GL_{n}(\Q_p)}W(g)W'(\delta_ng)\Phi(e_ng)|\det g|_p^{s}\,dg.
\end{align*}
Here $\delta_{n'} = {\rm diag}(-1,1,\cdots,(-1)^{n'})$, $e_{n} = (0,\cdots,0,1)$, and $dg$ is the quotient of the Haar measures on $\GL_{n'}(\Q_p)$ and $N_{n'}(\Q_p)$ with ${\rm vol}(\GL_{n'}(\Z_p)) = {\rm vol}(N_{n'}(\Z_p))=1$.
The integrals converge absolutely for ${\rm Re}(s)$ sufficiently large and admit meromorphic continuation to the whole complex plane. Moreover, they are represented by rational functions in $\C(p^{-s})$. The $\C$-vector space spanned by the local zeta integrals is a fractional $\C[p^s,p^{-s}]$-ideal of $\C(p^{-s})$ containing $1$. The $L$-factor $L(s,\itPi_p\times\itSigma_p)$ is then a generator of this fractional ideal.
We have the following functional equation:
If $n>n'$, we have
\begin{align}\label{E:fe 1}
\frac{Z^\vee(1-s,\rho(w_{n,n'})W^\vee,(W')^\vee)}{L(1-s,\itPi_p^\vee \times \itSigma_p^\vee)} = \omega_{\itSigma_p}(-1)^{n-1}\cdot\varepsilon(s,\itPi_p \times \itSigma_p,\psi_p)\cdot\frac{Z(s,W,W')}{L(s,\itPi_p \times \itSigma_p)}.
\end{align}
If $n=n'$, we have
\begin{align}\label{E:fe 2}
\frac{Z^\vee(1-s,W^\vee,(W')^\vee,\widehat{\Phi})}{L(1-s,\itPi_p^\vee \times \itSigma_p^\vee)} = \omega_{\itSigma_p}(-1)^{n-1}\cdot\varepsilon(s,\itPi_p \times \itSigma_p,\psi_p)\cdot\frac{Z(s,W,W',\Phi)}{L(s,\itPi_p \times \itSigma_p)}.
\end{align}
Here $W^\vee(g) = W(w_{n}{}^t\!g^{-1})$, $(W')^\vee(g) = W'(w_{n'}{}^t\!g^{-1})$, $w_{n,n'} = \bp {\bf 1}_{n'} & 0 \\ 0 & w_{n-n'}\ep$ with $w_N$ equal to the $N$ by $N$ anti-diagonal permutation matrix, and $\widehat{\Phi}$ is the Fourier transform given by
\[
\widehat{\Phi}(x) = \int_{\Q_p^n}\Phi(y) \psi_p(y{}^t\!x)\,dy.
\]
Let $\sigma \in {\rm Aut}(\C)$. Denote by ${}^\sigma\!\itPi_p$ the $\sigma$-conjugate of $\itPi_p$ and let $t_\sigma : \mathcal{W}(\itPi_p) \rightarrow \mathcal{W}({}^\sigma\!\itPi_p)$ be the $\sigma$-linear isomorphism defined by
\[
t_\sigma W(g) = \sigma \left( W({\rm diag}(u_{\sigma,p}^{1-n},u_{\sigma,p}^{2-n},\cdots,1)g) \right),\quad g \in \GL_n(\Q_p),
\]
where $u_{\sigma,p} \in \Z_p^\times$ is the element such that $\sigma (\psi_p(x)) = \psi_p(u_{\sigma,p}x)$ for $x \in \Q_p$.
We define $t_\sigma : \mathcal{W}(\itSigma_p) \rightarrow \mathcal{W}({}^\sigma\!\itSigma_p)$ is a similar way.
For a rational function $P \in \C(X)$, let ${}^\sigma\!P$ be the rational function obtained by acting $\sigma$ on the coefficients of $P$.
The following lemma is on the Galois equivariance properties of the local factors. We generalize the result of Waldspurger \cite[Proposition I.2.5]{Wald1985B} for $\GL_2 \times \GL_1$. For the local $L$-factors, one can also prove the equivariance following the arguments of Raghuram in \cite[Proposition 3.17]{Raghuram2009}.

\begin{lemma}\label{L:local factors}
For $\sigma \in {\rm Aut}(\C)$, we have the following identities as rational functions in $p^{-s}$.
\begin{align*}
{}^\sigma\!L(s+\tfrac{n-n'}{2},\itPi_p \times \itSigma_p) &=L(s+\tfrac{n'-n}{2},{}^\sigma\!\itPi_p \times {}^\sigma\!\itSigma_p),\\
{}^\sigma\!\varepsilon(s+\tfrac{n-n'}{2},\itPi_p \times \itSigma_p,\psi_p) &= 
{}^\sigma\!\omega_{\itPi_p}(u_{\sigma,p})^{n'}\cdot{}^\sigma\!\omega_{\itSigma_p}(u_{\sigma,p})^{n} \cdot\varepsilon(s+\tfrac{n-n'}{2},{}^\sigma\!\itPi_p \times {}^\sigma\!\itSigma_p,\psi_p).
\end{align*}
Here $u_{\sigma,p} \in \Z_p^\times$ is the element such that $\sigma (\psi_p(x)) = \psi_p(u_{\sigma,p}x)$ for $x \in \Q_p$.
\end{lemma}

\begin{proof}
We drop the subscript $p$ for brevity. 
Fix $\sigma \in {\rm Aut}(\C)$.
Let $W \in {\mathcal{W}(\itPi)}$, $W' \in {\mathcal{W}(\itSigma)}$, and $\Phi \in {S}(\Q_p^n)$. First we show that 
\begin{align}
{}^\sigma\!Z(s+\tfrac{n-n'}{2},W,W') &= Z(s+\tfrac{n-n'}{2},{}^\sigma W, {}^\sigma W'),\label{E:local factor pf 1}\\
{}^\sigma\!Z^\vee(s+\tfrac{n-n'}{2},W,W') &= Z^\vee(s+\tfrac{n-n'}{2},{}^\sigma W, {}^\sigma W')\label{E:local factor pf 2}
\end{align}
if $n > n'$, and 
\begin{align}
{}^\sigma\! Z(s,W,W',\Phi) = Z(s,{}^\sigma W,{}^\sigma W',{}^\sigma\!\Phi)\label{E:local factor pf 3}
\end{align}
if $n=n'$. 
Here ${}^\sigma W(g) = \sigma ( W(g))$, ${}^\sigma W'(g') = \sigma ( W'(g'))$, and ${}^\sigma\!\Phi(x) = \sigma(\Phi(x))$ for $g \in \GL_n(\Q_p)$, $g' \in \GL_{n'}(\Q_p)$, and $x \in \Q_p^n$.
By \cite[Proposition 2.2]{JPSS1979}, there exist a finite set $\frak{X}_\itPi$ of characters of $T_n(\Q_p)$ and a non-negative integer $N_{\itPi}$ such that
\[
W(ak) = \sum_{\underline{m} = (m_1,\cdots,m_n) \atop 0 \leq m_i \leq N_\itPi}\sum_{\chi \in \frak{X}_\itPi} \chi(a)\prod_{i=1}^n(\log_p|a_i|)^{-m_i}\cdot\Phi_{\underline{m},\chi}(a_1,\cdots,a_n;k)
\]
for all $a = {\rm diag}(a_1\cdots a_n,\,a_2\cdots a_n,\,...,\,a_n) \in T_n(\Q_p)$ and $k \in \GL_n(\Z_p)$, where $\Phi_{\underline{m},\chi} \in S(\Q_p^n\times \GL_n(\Z_p))$ is some Bruhat--Schwartz function depending on $\underline{m}$, $\chi$, and $W$.
Similarly, we have
\[
W'(\delta_{n'}a'k') = \sum_{\underline{m}' = (m_1',\cdots,m_{n'}') \atop 0 \leq m_i' \leq N_\itSigma}\sum_{\chi' \in \frak{X}_\itSigma} \chi'(a')\prod_{i=1}^{n'}(\log_p|a_i'|)^{-m_i'}\cdot\Phi_{\underline{m}',\chi'}(a_1',\cdots,a_{n'}';k')
\]
for all $a' = {\rm diag}(a_1'\cdots a_{n'}',\,a_2'\cdots a_{n'}',\,...,\,a_{n'}') \in T_{n'}(\Q_p)$ and $k' \in \GL_{n'}(\Z_p)$.
Therefore, when $n>n'$, we have
\begin{align*}
&Z(s+\tfrac{n-n'}{2},W,W')\\
&= \int_{\GL_{n'}(\Z_p)}dk'\,\prod_{1 \leq i \leq n'}\int_{\Q_p^\times}d^\times a_i'\, W\left(\bp a'k' & 0 \\ 0 & {\bf 1}_{n-n'}\ep\right)W'(\delta_{n'}a'k')\prod_{i=1}^{n'}|a_i'|^{i(s+i-m)}\\
& = \sum_{\underline{m},\underline{m}',\chi,\chi'}\int_{\GL_{n'}(\Z_p)}dk'\,Z(s,\underline{m},\underline{m}',\chi,\chi',\Phi_{\underline{m},\chi},\Phi_{\underline{m}',\chi'};k'),
\end{align*}
where 
\begin{align*}
&Z(s,\underline{m},\underline{m}',\chi,\chi',\Phi_{\underline{m},\chi},\Phi_{\underline{m}',\chi'};k')\\
& = \prod_{1 \leq i \leq n'}\int_{\Q_p^\times}d^\times a_i'\, \chi\left(\bp a' & 0 \\ 0 & {\bf 1}_{n-n'}\ep\right)\chi'(a')\prod_{i=1}^{n'}(\log_p|a_i'|)^{-m_i-m_i'}\\
&\quad\quad\quad\quad\quad\quad\Phi_{\underline{m},\chi}\left(a_1',\cdots,a_{n'}',1,\cdots,1;\bp k' & 0 \\ 0 & {\bf 1}_{n-n'} \ep \right)\Phi_{\underline{m}',\chi'}(a_1',\cdots,a_{n'}';k')\prod_{i=1}^{n'}|a_i'|^{i(s+i-m)}.
\end{align*}
Note that the factor $\prod_{i=1}^{n'}|a_i'|^{i(i-m)}$ is the Jacobian appearing in the integration formula with respect to the Iwasawa decomposition.
By the definition of ${}^\sigma W$, we have
\[
{}^\sigma W(ak) = \sum_{\underline{m} = (m_1,\cdots,m_n) \atop 0 \leq m_i \leq N_\itPi}\sum_{\chi \in \frak{X}_\itPi} {}^\sigma\!\chi(a)\prod_{i=1}^n(\log_p|a_i|)^{-m_i}\cdot{}^\sigma\!\Phi_{\underline{m},\chi}(a_1,\cdots,a_n;k)
\]
for all $a = {\rm diag}(a_1\cdots a_n,a_2\cdots a_n,\cdots,a_n) \in T_n(\Q_p)$ and $k \in \GL_n(\Z_p)$.
Similar formula holds for ${}^\sigma W'$.
Hence we also have
\[
Z(s+\tfrac{n-n'}{2},{}^\sigma W, {}^\sigma W') = \sum_{\underline{m},\underline{m}',\chi,\chi'}\int_{\GL_{n'}(\Z_p)}dk'\,Z(s,\underline{m},\underline{m}',{}^\sigma\!\chi,{}^\sigma\!\chi',{}^\sigma\!\Phi_{\underline{m},\chi},{}^\sigma\!\Phi_{\underline{m}',\chi'};k').
\]
On the other hand, $Z(s,\underline{m},\underline{m}',\chi,\chi',\Phi_{\underline{m},\chi},\Phi_{\underline{m}',\chi'};k')$ is a (generalized) Tate integral which represents a rational function in $p^{-s}$ whose denominator depends only on $\underline{m},\underline{m}',\chi,\chi'$.
Moreover, it is easy to verify that (cf.\,\cite[Proposition A]{Grobner2018})
\[
{}^\sigma\!Z(s,\underline{m},\underline{m}',\chi,\chi',\Phi_{\underline{m},\chi},\Phi_{\underline{m}',\chi'};k') = Z(s,\underline{m},\underline{m}',{}^\sigma\!\chi,{}^\sigma\!\chi',{}^\sigma\!\Phi_{\underline{m},\chi},{}^\sigma\!\Phi_{\underline{m}',\chi'};k')
\]
as rational functions in $p^{-s}$.
We thus obtain (\ref{E:local factor pf 1}). 
By \cite[Lemma 2.6]{JPSS1983}, 
\[
W\left( \bp g & 0 & 0 \\x & {\bf 1}_{n-n'-1} & 0 \\ 0 & 0 & 1\ep \right)\neq 0
\]
implies that $x$ belongs to a compact set of $M_{n-n'-1,n'}(\Q_p)$ which is independent of $g \in \GL_{n'}(\Q_p)$. 
Similarly we can show that (\ref{E:local factor pf 2}) and (\ref{E:local factor pf 3}) hold. 
The Galois equivariance of local $L$-factors then follows immediately from (\ref{E:local factor pf 1}) and (\ref{E:local factor pf 3}).
For the equivariance of local $\varepsilon$-factors, we begin with the following equalities: If $n>n'$ and $W$ is right invariant by 
\[
\mathcal{T}_1 = \bp {\bf 1}_{n'} & 0 & 0 \\ 0 & T_{n-n'-1}(\Z_p) & 0 \\ 0 & 0 & 1 \ep,
\]
then
\begin{align}\label{E:local factor pf 4}
\begin{split}
{}^\sigma\!Z(s+\tfrac{n-n'}{2},W,W') &= {}^\sigma\!\omega_{\itSigma}(u_\sigma)^{n'-n}\cdot Z(s+\tfrac{n-n'}{2},t_\sigma W,t_\sigma W'),\\
{}^\sigma\!Z^\vee(s+\tfrac{n-n'}{2},W,W') &= {}^\sigma\!\omega_{\itSigma}(u_\sigma)^{n'-n}\cdot Z^\vee(s+\tfrac{n-n'}{2},t_\sigma W,t_\sigma W').
\end{split}
\end{align}
If $n>n'$ and $W$ is right invariant by 
\[
\mathcal{T}_2 = \bp {\bf 1}_{n'+1} & 0  \\ 0 & T_{n-n'-1}(\Z_p)  \ep,
\]
then
\begin{align}\label{E:local factor pf 5}
\begin{split}
{}^\sigma\!Z(s+\tfrac{n-n'}{2},W,W') &= {}^\sigma\!\omega_{\itPi}(u_\sigma)^{n-n'-1}\cdot {}^\sigma\!\omega_\itSigma(u_\sigma)^{-1}\cdot Z(s+\tfrac{n-n'}{2},t_\sigma W,t_\sigma W'),\\
{}^\sigma\!Z^\vee(s+\tfrac{n-n'}{2},W,W') &= {}^\sigma\!\omega_{\itPi}(u_\sigma)^{n-n'-1}\cdot{}^\sigma\!\omega_\itSigma(u_\sigma)^{-1}\cdot Z^\vee(s+\tfrac{n-n'}{2},t_\sigma W,t_\sigma W').
\end{split}
\end{align}
If $n=n'$, then
\begin{align}\label{E:local factor pf 6}
{}^\sigma\!Z(s,W,W',\Phi) = Z(s,t_\sigma W,t_\sigma W',{}^\sigma\!\Phi).
\end{align}
By (\ref{E:local factor pf 1}) and (\ref{E:local factor pf 2}), (\ref{E:local factor pf 4}) follows from the change of variables $g \mapsto {\rm diag}(u_\sigma^{1-n},u_{\sigma}^{2-n}\cdots,u_\sigma^{n'-n})\cdot g$, the $\mathcal{T}_1$-invariance of $W$, and 
\[
{}^\sigma W'({\rm diag}(u_\sigma^{1-n},u_{\sigma}^{2-n}\cdots,u_\sigma^{n'-n})\cdot g) = {}^\sigma\!\omega_{\itSigma}(u_\sigma)^{n'-n}\cdot t_\sigma W'(g).
\]
Similarly, (\ref{E:local factor pf 5}) follows from the change of variables $g \mapsto {\rm diag}(u_\sigma^{-n'},u_{\sigma}^{1-n'}\cdots,u_\sigma^{-1})\cdot g$ and 
\begin{align*}
{}^\sigma W'( {\rm diag}(u_\sigma^{-n'},u_{\sigma}^{1-n'}\cdots,u_\sigma^{-1})\cdot g) &= {}^\sigma\!\omega_{\itSigma}(u_\sigma)^{-1}\cdot t_\sigma W'(g),\\
{}^\sigma W\left( \bp {\rm diag}(u_\sigma^{-n'},u_{\sigma}^{1-n'}\cdots,u_\sigma^{-1})\cdot g & 0 \\ 0 & {\bf 1}_{n-n'} \ep\right) &= {}^\sigma\!\omega_{\itPi}(u_\sigma)^{n-n'-1}\cdot t_\sigma W\left( \bp g & 0 \\ 0 &{\bf 1}_{n-n'}\ep\right).
\end{align*}
Here the second equality follows from the $\mathcal{T}_2$-invariance of $W$. Also (\ref{E:local factor pf 6}) follows from (\ref{E:local factor pf 3}) and the change of variables $g \mapsto {\rm diag}(u_\sigma^{1-n'},u_{\sigma}^{2-n'}\cdots,1)\cdot g$.
By definition, we easily see that
\begin{align*}
t_\sigma W^\vee = {}^\sigma\!\omega_\itPi(u_\sigma)^{n-1}\cdot (t_\sigma W)^\vee,\quad
t_\sigma(W')^\vee = {}^\sigma\!\omega_\itSigma(u_\sigma)^{n'-1}\cdot (t_\sigma W')^\vee.
\end{align*}
Therefore, if $n>n'$ and $W$ is right $\mathcal{T}_1$-invariant (hence $\rho(w_{n,n'})W^\vee$ is right $\mathcal{T}_2$-invariant), by (\ref{E:local factor pf 5}) we have 
\begin{align}\label{E:local factor pf 7}
\begin{split}
&{}^\sigma\!Z^\vee(1-s-\tfrac{n-n'}{2},\rho(w_{n,n'})W^\vee,(W')^\vee)\\
&= {}^\sigma\!\omega_{\itPi}(u_\sigma)^{-n+n'+1}\cdot {}^\sigma\!\omega_\itSigma(u_\sigma)\cdot Z(s+\tfrac{n-n'}{2},\rho(w_{n,n'})t_\sigma W^\vee,t_\sigma (W')^\vee)\\
&= {}^\sigma\!\omega_{\itPi}(u_\sigma)^{n'}\cdot {}^\sigma\!\omega_\itSigma(u_\sigma)^{n'}\cdot Z(s+\tfrac{n-n'}{2},\rho(w_{n,n'})(t_\sigma W)^\vee,(t_\sigma W')^\vee).
\end{split}
\end{align}
Similarly, if $n=n'$, by (\ref{E:local factor pf 6}) and ${}^\sigma\!(\widehat{\Phi})(x) = \widehat{{}^\sigma\!\Phi}(u_\sigma x)$ we have
\begin{align}\label{E:local factor pf 8}
\begin{split}
&{}^\sigma\!Z(1-s,W^\vee,(W')^\vee,\widehat{\Phi}) = {}^\sigma\!\omega_{\itPi}(u_\sigma)^{n}\cdot {}^\sigma\!\omega_\itSigma(u_\sigma)^n\cdot Z(1-s,(t_\sigma W)^\vee,(t_\sigma W')^\vee,\widehat{{}^\sigma\!\Phi}).
\end{split}
\end{align}
In this case, the Galois equivariance of local $\varepsilon$-factors follows from (\ref{E:local factor pf 6}) and (\ref{E:local factor pf 8}) by applying $\sigma$ to both sides of the functional equation (\ref{E:fe 2}).
If $n>n'$, then the equivariance is a consequence of the functional equation (\ref{E:fe 1}) and (\ref{E:local factor pf 4}), (\ref{E:local factor pf 7}) provided we can show that there exist $\mathcal{T}_1$-invariant $W$ and some $W'$ such that $Z(s,W,W')$ is non-zero.  
The existence of such $W$ and $W'$ is a simple application of the well-known result of Gelfand and Kazhdan \cite[Proposition 2]{GK1975}.
Indeed, let $P_{n}$ be the mirabolic subgroup of $\GL_{n}$ consisting of matrices whose last rows are $e_{n}$.
Fix $W'$ such that $W'(\delta_{n'}) =1$ and let $m$ be a sufficiently large integer such that $W'$ is right invariant by ${\bf 1}_{n'}+p^mM_{n',n'}(\Z_p)$.
By $loc.$ $cit.$, there exists a Whittaker function $W$ uniquely determined so that 
\[
W\vert_{P_{n}(\Q_p)}(g) = \begin{cases}
0 & \mbox{ if $g \notin N_{n}(\Q_p) \bp K_m & 0 \\ 0 & 1  \ep$},\\
\psi_{N_{n}}(u) & \mbox{ if $g = u\cdot k \in N_{n}(\Q_p) \bp K_m & 0 \\ 0 & 1  \ep$},
\end{cases}
\]
where $K_m$ is an open compact subgroup of $\GL_{n-1}(\Q_p)$ given by
\[
K_m = \bp {\bf 1}_{n'}+ p^mM_{n',n'}(\Z_p) & p^m M_{n',n-n'-1}(\Z_p) \\  p^m M_{n-n'-1,n'}(\Z_p) & T_{n-n'-1}(\Z_p) + p^mM_{n-n'-1,n-n'-1}(\Z_p) \ep.
\]
Then we have
\[
Z(s,W,W') = W'(\delta_{n'})\cdot \int_{N_{n'}(\Q_p)\backslash \GL_{n'}(\Q_p)}W\left( \bp g & 0 \\ 0 & {\bf 1}_{n-n'} \ep \right)\,dg \neq 0.
\]
This completes the proof.
\end{proof}




\subsection{Ratios of critical values under duality}

Following is the main result of this section. Similar result was proved in \cite[Appendix A]{Chen2021} for standard $L$-functions.

\begin{thm}\label{T:ratio under duality}
Let $m_0 \in \Z+\tfrac{n+n'}{2}$ be a critical point for $L(s,\itPi \times \itSigma)$ such that $L(m_0,{}^\sigma\!\itPi \times {}^\sigma\!\itSigma) \neq 0$ for all $\sigma \in {\rm Aut}(\C)$.  
Then we have
\begin{align*}
&\sigma \left( \frac{L(m_0,\itPi \times \itSigma)}{\varepsilon(0,\itPi_\infty\times\itSigma_\infty,\psi_\infty)\cdot G(\omega_\itPi)^{n'}\cdot G(\omega_\itSigma)^{n}\cdot L(1-m_0,\itPi^\vee \times \itSigma^\vee)} \right)\\
& = \frac{L(m_0,{}^\sigma\!\itPi \times {}^\sigma\!\itSigma)}{\varepsilon(0,\itPi_\infty\times\itSigma_\infty,\psi_\infty)\cdot G({}^\sigma\!\omega_\itPi)^{n'}\cdot G({}^\sigma\!\omega_\itSigma)^n\cdot L(1-m_0,{}^\sigma\!\itPi^\vee \times {}^\sigma\!\itSigma^\vee)},\quad \sigma \in {\rm Aut}(\C).
\end{align*}
\end{thm}

\begin{proof}
By the global functional equation (\ref{E:global fe}), we have
\begin{align*}
\frac{L(m_0,\itPi \times \itSigma)}{\varepsilon(m_0,\itPi_\infty \times \itSigma_\infty,\psi_\infty)\cdot L(1-m_0,\itPi^\vee \times \itSigma^\vee)} = \prod_p\varepsilon(m_0,\itPi_p \times \itSigma_p,\psi_p).
\end{align*}
Let $\sigma \in {\rm Aut}(\C)$.
Since $\varepsilon(s,\itPi_p \times \itSigma_p,\psi_p) = 1$ when $\itPi_p$ and $\itSigma_p$ are unramified, together with (\ref{E:Galois Gauss sum}) and Lemma \ref{L:local factors}, we have
\begin{align*}\label{E:main 1 proof 1}
&\sigma \left(\frac{L(m_0,\itPi \times \itSigma)}{\varepsilon(m_0,\itPi_\infty \times \itSigma_\infty,\psi_\infty)\cdot L(1-m_0,\itPi^\vee \times \itSigma^\vee)}\right)\\
& = \prod_p\sigma(\varepsilon(m_0,\itPi_p \times \itSigma_p,\psi_p))\\
& = \prod_p
{}^\sigma\!\omega_{\itPi_p}(u_{\sigma,p})^{n'} \cdot {}^\sigma\!\omega_{\itSigma_p}(u_{\sigma,p})^{n}\cdot\varepsilon(m_0,{}^\sigma\!\itPi_p \times {}^\sigma\!\itSigma_p,\psi_p)\\
& = {}^\sigma\!\omega_\itPi(u_\sigma)^{n'}\cdot {}^\sigma\!\omega_\itSigma(u_\sigma)^n\cdot \prod_p \varepsilon(m_0,{}^\sigma\!\itPi_p \times {}^\sigma\!\itSigma_p,\psi_p)\\
& = \frac{\sigma (G(\omega_\itPi)^{n'}\cdot G(\omega_\itSigma)^n)}{G({}^\sigma\!\omega_\itPi)^{n'}\cdot G({}^\sigma\!\omega_\itSigma)^n}\cdot\prod_p \varepsilon(m_0,{}^\sigma\!\itPi_p \times {}^\sigma\!\itSigma_p,\psi_p).
\end{align*}
By the global functional equation (\ref{E:global fe}) for $L(s,{}^\sigma\!\itSigma \times {}^\sigma\!\itPi)$, we have
\begin{align*}
\frac{L(m_0,{}^\sigma\!\itPi \times {}^\sigma\!\itSigma)}{\varepsilon(m_0,\itPi_\infty \times \itSigma_\infty,\psi_\infty)\cdot L(1-m_0,{}^\sigma\!\itPi^\vee \times {}^\sigma\!\itSigma^\vee)} = \prod_p\varepsilon(m_0,{}^\sigma\!\itPi_p \times {}^\sigma\!\itSigma_p,\psi_p).
\end{align*}
We thus conclude that
\begin{align*}
&\sigma \left( \frac{L(m_0,\itPi \times \itSigma)}{\varepsilon(m_0,\itPi_\infty\times\itSigma_\infty,\psi_\infty)\cdot G(\omega_\itPi)^{n'}\cdot G(\omega_\itSigma)^n\cdot L(1-m_0,\itPi^\vee \times \itSigma^\vee)} \right)\\
& = \frac{L(m_0,{}^\sigma\!\itPi \times {}^\sigma\!\itSigma)}{\varepsilon(m_0,\itPi_\infty\times\itSigma_\infty,\psi_\infty)\cdot G({}^\sigma\!\omega_\itPi)^{n'}\cdot G({}^\sigma\!\omega_\itSigma)^n\cdot L(1-m_0,{}^\sigma\!\itPi^\vee \times {}^\sigma\!\itSigma^\vee)}.
\end{align*}
Finally, by our choice of $\psi_\infty$, we have
$\varepsilon(s,\itPi_\infty\times\itSigma_\infty,\psi_\infty) = \varepsilon(0,\itPi_\infty\times\itSigma_\infty,\psi_\infty)$ (cf.\,\cite[\S\,16]{Jacquet2009}).
This completes the proof.
\end{proof}

\section{Proof of main results}\label{S:main results}

The aim of this section is to prove our main results Theorems \ref{T:main 1} and \ref{T:main 2}. 

\subsection{Essentially self-dual cuspidal automorphic representations}\label{SS:polarized}

Let $\itPi$ be a cuspidal automorphic representation of $\GL_N(\A)$. We say $\itPi$ is essentially self-dual if there exists a Hecke character $\chi$ such that $\itPi = \itPi^\vee \otimes \chi$. In this case, we have
\[
L(s,\itPi \times \itPi^\vee) = L(s, \itPi, {\rm Sym}^2 \otimes \chi^{-1})\cdot L(s, \itPi, \text{\large$ \wedge $}^2 \otimes \chi^{-1}),
\]
and exactly one of the $L$-functions on the right-hand side has a pole at $s=1$.
We say $\itPi$ is $\chi$-orthogonal (resp.\,$\chi$-symplectic) if $L(s, \itPi, {\rm Sym}^2 \otimes \chi^{-1})$ (resp.\,$L(s, \itPi, \text{\large$ \wedge $}^2 \otimes \chi^{-1})$) has a pole at $s=1$.
If $N=2n+1$, then $\itPi$ must be $\chi$-orthogonal by the result of Jacquet and Shalika \cite[Theorem 2]{JS1990b}.
Also, $\itPi \otimes \chi^n\omega_\itPi^{-1}$ is self-dual with trivial central character and is the functorial transfer of a globally generic discrete automorphic representation of $\Sp_{2n}(\A)$ by the endoscopic classification of Arthur \cite{Arthur2013}. Moreover, the descent to $\Sp_{2n}(\A)$ is cuspidal by the descent method of Ginzburg, Rallis, and Soudry \cite{GRS2011}.
If $N=2n$, then $\chi^{n}\omega_\itPi^{-1}$ is a quadratic Hecke character. We consider the quasi-split general spin group ${\rm GSpin}_{2n}^*$ over $\Q$ determined by $\chi^{n}\omega_\itPi^{-1}$ as in \cite[\S\,2.1.3]{AS2014} and the split general spin group ${\rm GSpin}_{2n+1}$.
By the results of Asgari--Shahidi \cite{AS2006b},\cite{AS2014} and Hundley--Sayag \cite{HS2016}, $\itPi$ is the weak functorial transfer of a globally generic cuspidal automorphic representation of ${\rm GSpin}_{2n}^*(\A)$ (resp.\,${\rm GSpin}_{2n+1}(\A)$) with central character $\chi$ if $\itPi$ is $\chi$-orthogonal (resp.\,$\chi$-symplectic).
Moreover, the transfer is strong at the archimedean place.
In the following lemma, when $\itPi$ is regular algebraic and $N=2n$, we show that the type of $\itPi$ with respect $\chi$ is determined by $\varepsilon(\chi_\infty)$.

\begin{lemma}\label{L:auxiliary 3}
Let $\itPi$ be a regular algebraic cuspidal automorphic representation of $\GL_{2n}(\A)$ with $\itPi = \itPi^\vee \otimes \chi$. Then $\itPi$ is $\chi$-orthogonal if and only if $\varepsilon(\chi_\infty)=-1$.
\end{lemma}

\begin{proof}
Let $(\underline{\kappa};\,{\sf w})$ be the infinity type of $\itPi_\infty$. Then the local $L$-parameter $\phi_{\itPi_\infty}$ of $\itPi_\infty$ is given by
\[
\phi_{\itPi_\infty} = \phi_{\kappa_1} \otimes |\mbox{ }|^{{\sf w}/2} \oplus \cdots \oplus\phi_{\kappa_n} \otimes |\mbox{ }|^{{\sf w}/2}. 
\]
Here $\phi_\kappa$ is the irreducible $2$-dimensional representation of $W_\R$ corresponding to the discrete series representation of $\GL_2(\R)$ of minimal weight $\kappa \geq 2$ via the local Langlands correspondence.
If $\itPi$ is $\chi$-orthogonal,
by the results of Asgari--Shahidi and Hundley--Sayag recalled above, the image of $\phi_{\itPi_\infty}$ factors through ${\rm GO}_{2n}(\C) \subset \GL_{2n}(\C)$ and 
\[
{\rm Hom}_{W_\R}({\rm Sym}^2\phi_{\itPi_\infty},\chi_\infty) \neq 0.
\]
For $i \neq j$, we have $\kappa_i \neq \kappa_j$ and 
\begin{align}\label{E:local tensor}
\phi_{\kappa_i} \otimes \phi_{\kappa_j} = \phi_{\kappa_i+\kappa_j-1}\oplus \phi_{|\kappa_i-\kappa_j|+1}.
\end{align}
Also
\[
{\rm Sym}^2\phi_{\kappa_i} = \phi_{2\kappa_i-1} \oplus {\rm sgn}^{\kappa_i-1}.
\]
Since $\kappa_i \equiv {\sf w}\,({\rm mod}\,2)$, we see that
\[
{\rm Hom}_{W_\R}({\rm Sym}^2\phi_{\itPi_\infty},\chi_\infty) = \bigoplus_{i=1}^n {\rm Hom}_{W_\R}({\rm sgn}^{{\sf w}-1}|\mbox{ }|^{\sf w},\chi_\infty).
\]
Thus we have $\chi_\infty = {\rm sgn}^{{\sf w}-1}|\mbox{ }|^{\sf w}$ which implies that $\varepsilon(\chi_\infty)=-1$.
If $\itPi$ is $\chi$-symplectic, then the image of $\phi_{\itPi_\infty}$ factors through $\GSp_{2n}(\C) \subset \GL_{2n}(\C)$ and 
\[
{\rm Hom}_{W_\R}(\text{\large$ \wedge $}^2\phi_{\itPi_\infty},\chi_\infty) = \bigoplus_{i=1}^n {\rm Hom}_{W_\R}({\rm sgn}^{{\sf w}}|\mbox{ }|^{\sf w},\chi_\infty) \neq 0.
\]
Hence $\chi_\infty = {\rm sgn}^{{\sf w}}|\mbox{ }|^{\sf w}$ and $\varepsilon(\chi_\infty)=1$.
This completes the proof.
\end{proof}

An isobaric automorphic representation
\[
\itPi = \itPi_1 \boxplus \cdots \boxplus \itPi_k
\]
of $\GL_n(\A)$ is called $\chi$-orthogonal (resp.\,$\chi$-symplectic) if $\itPi_i \neq \itPi_j$ for $i \neq j$ and $\itPi_i$ is $\chi$-orthogonal (resp.\,$\chi$-symplectic) for all $i$.
Lemma \ref{L:auxiliary 3} gives a simple criterion for testing the types of regular algebraic isobaric automorphic representations under Langlands functoriality. For instance, we consider the automorphic tensor product and the Asai transfer.
For cuspidal automorphic representations $\itPi_1$ and $\itPi_2$ of $\GL_{n_1}(\A)$ and $\GL_{n_2}(\A)$, we denote by $\itPi_1 \boxtimes \itPi_2$ the automorphic tensor product of $\itPi_1$ and $\itPi_2$. Recall $\itPi_1 \boxtimes \itPi_2$ is the irreducible admissible $((\frak{g}_{n_1n_2},{\rm O}_{n_1n_2}(\R))\times \GL_{n_1n_2}(\A_f))$-module defined by the restricted tensor product
\[
\itPi_1 \boxtimes \itPi_2 = \bigotimes_v \itPi_{1,v}\boxtimes \itPi_{2,v}.
\]
Here $\itPi_{1,v}\boxtimes \itPi_{2,v}$ correspondes to $\phi_{\itPi_{1,v}} \otimes \phi_{\itPi_{2,v}}$ via the local Langlands correspondence and $\phi_{\itPi_{i,v}}$ is the $L$-parameter of $\itPi_{i,v}$ for $i=1,2$.
The functoriality of the automorphic tensor products is known for $n_1=n_2=2$ \cite{Rama2000} or $n_1=2$, $n_2=3$ \cite{KS2002} (see also \cite{Dieulefait2020}, \cite{AdD2022}).

\begin{corollary}\label{C:tensor}
Let $\itPi_1$ and $\itPi_2$ be regular algebraic cuspidal automorphic representations of $\GL_{n_1}(\A)$ and $\GL_{n_2}(\A)$ with $\itPi_1 = \itPi_1^\vee \otimes\chi_1$ and $\itPi_2 = \itPi_2^\vee \otimes\chi_2$.
Assume $n_1n_2$ is even, $\itPi_1 \boxtimes \itPi_2$ is regular and automorphic.
Then $(\itPi_1 \boxtimes \itPi_2) \otimes |\mbox{ }|_\A^{(n_1+n_2-1)/2}$ is regular algebraic and isobaric. Moreover, it is $\chi_1\chi_2|\mbox{ }|_\A^{n_1+n_2-1}$-orthogonal if and only if $\varepsilon(\chi_{1,\infty}\chi_{2,\infty}) = (-1)^{n_1+n_2}$.
\end{corollary}

\begin{proof}
The algebraicity of $(\itPi_1 \boxtimes \itPi_2) \otimes |\mbox{ }|_\A^{(n_1+n_2-1)/2}$ follows immediately from the assumption that $n_1n_2$ is even and (\ref{E:local tensor}).
Since $\itPi_1$ and $\itPi_2$ are regular algebraic and essentially self-dual, they are essentially tempered everywhere. 
Indeed, by base change to imaginary quadratic extensions and \cite[Lemma 4.1.4]{CHT2008}, we are reduced to the temperedness of regular algebraic conjugate self-dual cuspidal automorphic representations over CM-fields.
We refer to \cite[Theorem 1.2]{Caraiani2012} and the references therein.
Hence $\itPi_1\boxtimes\itPi_2$ is also essentially tempered and 
\[
\itPi_1\boxtimes\itPi_2 = \itSigma_1 \boxplus \cdots \boxplus \itSigma_k
\]
for some cuspidal automorphic representation $\itSigma_i$ of $\GL_{N_i}(\A)$ (cf.\,\cite[Lemme 1.5]{Clozel1990}).
Moreover, the proof in $loc.$ $cit.$ shows that there exists $t \in \C$ such that $\itSigma_i\otimes|\mbox{ }|_\A^t$ is unitary for all $i$.
The regularity assumption on $\itPi_1 \boxtimes \itPi_2$ then implies that $\itSigma_i \neq \itSigma_j$ for $i \neq j$. Also note that $N_i$ is even for all $i$. Indeed, if $N_i$ is odd for some $i$, then $N_j$ is also odd for some $i \neq j$ as $\sum_{k}N_k = n_1n_2$ is even. In this case, up to twisting by a sign character, the Langlands parameters of $\itSigma_{i,\infty}$ and $\itSigma_{j,\infty}$ have a common $1$-dimensional sub-representation of $W_\R$, which contradicts the regularity of $\itPi_1 \boxtimes \itPi_2$.
The essentially self-dual condition implies that
\[
(\itPi_1\boxtimes\itPi_2) = (\itPi_1\boxtimes\itPi_2)^\vee\otimes \chi_1\chi_2.
\]
The regularity assumption again implies that $\itSigma_i = \itSigma_i^\vee\otimes\chi_1\chi_2$ for all $i$.
Since $N_i$ is even and $\itSigma_i\otimes |\mbox{ }|_\A^{(n_1+n_2-1)/2}$ is regular algebraic for all $i$, the assertion then follows from Lemma \ref{L:auxiliary 3}.
This completes the proof.  
\end{proof}

Let $\F$ be a number field with $d=[\F:\Q]$. For each place $v$ of $\Q$, let $\F_v = \prod_{w\mid v}\F_w$ be the product of localizations of $\F$ at places dividing $v$, and $d_w = [\F_w:\Q_v]$ for $w \mid v$.
Let $\pi$ be a cuspidal automorphic representation of $\GL_N(\A_\F) = ({\rm R}_{\F/\Q}\GL_{N})(\A)$.
For each place $v$ of $\Q$, we have the associated $L$-parameter 
\[
\phi_{\pi_v} : W_{\Q_v}' \longrightarrow \prod_{w\mid v}\GL_N(\C)^{d_w}\rtimes \Gal(\overline{\Q}_v/\Q_v),
\]
where $W_{\Q_v}'$ is the Weil--Deligne group of $\Q_v$ and the action of $\Gal(\overline{\Q}_v/\Q_v)$ on $\prod_{w\mid v}\GL_N(\C)^{d_w}$ is the permutation of components induced by the natural homomorphism $\Gal(\overline{\Q}_v/\Q_v) \rightarrow \prod_{w \mid v}\Gal(\F_w^{\Gal}/\Q_v)$.
Let ${\rm As}_v$ be the Asai representation of $\prod_{w\mid v}\GL_N(\C)^{d_w}\rtimes \Gal(\overline{\Q}_v/\Q_v)$ on $\otimes_{w \mid v}(\C^N)^{\otimes d_w}$ defined so that
\[
{\rm As}_v(\prod_{w \mid v}g_w)\cdot (\otimes_{w \mid v}{\bf v}_w) = \otimes_{w \mid v}g_w\cdot {\bf v}_w,\quad g_w \in \GL_N(\C)^{d_w},\,{\bf v}_w \in (\C^N)^{\otimes d_w}
\]
and the action of $\Gal(\overline{\Q}_v/\Q_v)$ on $\otimes_{w \mid v}(\C^N)^{\otimes d_w}$ is the permutation of components induced by the natural homomorphism $\Gal(\overline{\Q}_v/\Q_v) \rightarrow \prod_{w \mid v}\Gal(\F_w^{\Gal}/\Q_v)$.
The Asai transfer ${\rm As}(\pi)$ of $\pi$ is the irreducible admissible $((\frak{g}_{N^d},{\rm O}_{N^d}(\R))\times \GL_{N^d}(\A_f))$-module defined by the restricted tensor product
\[
{\rm As}(\pi) = \bigotimes_{v} {\rm As}_v(\pi_v).
\]
Here ${\rm As}_v(\pi_v)$ correspondes to ${\rm As}_v\circ \phi_{\pi_v}$ via the local Langlands correspondence.
The functoriality of the Asai transfers is known for $d=N=2$ \cite{Kris2003}.

\begin{corollary}\label{C:Asai}
Let $\F$ be a totally real number field with $[\F:\Q]=d$.
Let $\pi$ be a regular algebraic cuspidal automorphic representation of $\GL_N(\A_\F)$ with $\pi = \pi^\vee \otimes \chi$. 
Assume $N$ is even, ${\rm As}(\pi)$ is regular and automorphic. Then ${\rm As}(\pi)\otimes |\mbox{ }|_\A^{(d-1)/2}$ is regular algebraic and isobaric. 
Moreover, it is $\chi\vert_{\A^\times}|\mbox{ }|_\A^{d-1}$-orthogonal if and only if $\varepsilon(\prod_{w \mid \infty}\chi_w) = (-1)^d$.
\end{corollary}

\begin{proof}
Since the archimedean place $\infty$ splits in $\F$, by definition we have
\[
{\rm As}_\infty(\pi_\infty) = \bigboxtimes_{w \mid \infty} \pi_w.
\]
The algebraicity of ${\rm As}(\pi)\otimes |\mbox{ }|_\A^{(d-1)/2}$ then follows from (\ref{E:local tensor}).
By the essentially self-dual condition and \cite[Lemma 7.1-(a)-(c)]{Prasad1992}, we have
\[
{\rm As}(\pi) = {\rm As}(\pi)^\vee \otimes \chi\vert_{\A^\times}.
\]
The rest of the proof is the same as the one for Corollary \ref{C:tensor}. 
\end{proof}

\begin{rmk}
If $\F$ is not totally real, then ${\rm As}(\pi)$ is not regular for any regular algebraic cuspidal automorphic representation $\pi$ of $\GL_N(\A_\F)$ with $N>1$.
\end{rmk}

\subsection{Auxiliary lemmas}

In this section, we recall two auxiliary lemmas which will be used in the proof of our main result Theorem \ref{T:main 1}.
We begin with the existence of regular algebraic self-dual cuspidal automorphic representations of $\GL_n(\A)$ proved by Bhagwat and Raghuram \cite[Theorem 2.10]{BR2017b}. 
In $loc.$ $cit.$, the cuspidality follows from \cite[Proposition 8.2]{MS2020} which is for $\Sp_{n-1}$ (hence $n$ is odd).
For the convenience of readers, we fill in the remaining details and give a uniform proof.

\begin{lemma}\label{L:auxiliary 1}
Let $(\underline{\kappa};\,0)$ be the infinity type of an irreducible admissible $(\frak{g}_n,{\rm O}_n(\R))$-module in $\Omega_0(n)$. Then there exists a regular algebraic self-dual cuspidal automorphic representation of $\GL_n(\A)$ with infinity type $(\underline{\kappa};\,0)$.
\end{lemma}

\begin{proof}
Let $r = \lfloor \tfrac{n}{2} \rfloor$, $G_r$ be the split odd special orthogonal group of rank $r$ if $n$ is even and $G_r = \Sp_{2r}$ if $n$ is odd.
With respect to the standard Cartan subgroup and positive system, let $\tau_{{\lambda}}$ be the discrete series representation of $G_r(\R)$ with Harish--Chandra parameter ${\lambda}$ with 
\[
{\lambda} = (\tfrac{\kappa_1-1}{2},\cdots,\tfrac{\kappa_r-1}{2}).
\]
It is proved in \cite{BR2017b} that the functorial transfer of $\tau_\lambda$ to $\GL_n(\R)$ belongs to $\Omega_0(n)$ with infinity type $(\underline{\kappa};\,0)$.
By the result of Clozel \cite[Theorem 1B]{Clozel1986}, there exists a cuspidal automorphic representation $\pi$ of $G_r(\A)$ such that 
\begin{itemize}
\item $\pi_p$ is the Steinberg representation of $G_r(\Q_p)$ for some prime $p$.
\item $\pi_\infty = \tau_\lambda$.
\end{itemize}
The global Arthur parameter of $\pi$ is of the form
\[
\itPsi = \itPi_1[d_1] \boxplus \cdots \boxplus \itPi_k[d_k]
\]
for some self-dual cuspidal automorphic representation $\itPi_i$ of $\GL_{n_i}(\A)$ and some $d_i \geq 1$ for $1 \leq i \leq k$ such that $\itPi_i$ is symplectic (resp.\,orthogonal) only if $d_i$ is odd when $n$ is even (resp.\,odd). 
By Arthur's multiplicity formula \cite[Theorems 1.5.1 and  1.5.2]{Arthur2013}, it suffices to show that 
\[
d_1=\cdots=d_k=1,\quad k=1. 
\]
Indeed, this would imply that $\itPi_1$ is the global functorial transfer of $\pi$ to $\GL_n(\A)$ whose archimedean component belongs to $\Omega_0(n)$ with infinity type $(\underline{\kappa};\,0)$.
Since $\pi_p$ is Steinberg, it is generic and its Langlands parameter $\phi_{\pi_p}$ is irreducible as a $n$-dimensional representation of the Weil--Deligne group of $\Q_p$.
By the result of Hazeltine, Liu, and Lo \cite[Theorem 1.8]{HLL2022} on enhanced Shahidi's conjecture, the genericity then implies that any local Arthur packet of $G_r(\Q_p)$ containing $\pi_p$ must be tempered.
Thus we have $d_1=\cdots=d_k=1$.
Therefore, for each place $v$ of $\Q$, the localization of $\itPsi$ at $v$ is equal to the local Langlands parameter of $\pi_v$. In particular, we have $\phi_{\pi_p} = \phi_{\itPi_{1,p}} \oplus \cdots \oplus \phi_{\itPi_{k,p}}$, where $\itPi_{i,p}$ refers to the local component of $\itPi_i$ at $p$.
We thus conclude from the irreducibility of $\phi_{\pi_p}$ that $k=1$.
This completes the proof.
\end{proof}

The following result is a relation between the archimedean periods defined in Definition \ref{D:archimedean periods}. 
The period relation was proved by Januszewski \cite{Januszewski2019} subject to the existence of rational cohomological test vectors.
In \cite{LLS2022}, Li, Liu, and Sun prove the period relation unconditionally. 
We give an independent proof based on global arguments.

\begin{lemma}\label{L:auxiliary 2}
Let $\itPi_\infty \in \Omega_0(n)$ and $\itSigma_\infty \in {\Omega_0(n-1)}$ with infinity types $(\underline{\kappa};\,{\sf w})$ and $(\underline{\ell};\,{\sf u})$ respectively. Assume $(\itPi_\infty,\itSigma_\infty)$ is balanced. Let $m_1,\,m_2 \in \Z$ such that $m_1 \neq \tfrac{-{\sf w}-{\sf u}}{2}$, $m_2 \neq \tfrac{-{\sf w}-{\sf u}}{2}$, and $L(s,\itPi_\infty \times \itSigma_\infty)$ and $L(1-s,\itPi_\infty^\vee \times \itSigma_\infty^\vee)$ are holomorphic at $s=m_1+\tfrac{1}{2},\,m_2+\tfrac{1}{2}$. Then we have
\begin{align}\label{E:local period}
\frac{p(m_1,\itPi_\infty \times \itSigma_\infty)}{p(m_2,\itPi_\infty \times \itSigma_\infty)} \in (\sqrt{-1})^{(m_1-m_2)n(n-1)/2}\cdot \Q^\times.
\end{align}
\end{lemma}

\begin{proof}
We assume $n$ is even. Thus ${\sf u}$ must be even.
The case when $n$ is odd can be proved in a similar way.
We prove the period relation by global arguments.
Assume there exist cuspidal automorphic representations $\itPi$ and $\itSigma$ of $\GL_{n}(\A)$ and $\GL_{n-1}(\A)$ with archimedean components $\itPi_\infty$ and $\itSigma_\infty$, respectively. The assumptions on $m_1$ and $m_2$ are equivalent to saying that $m_1+\tfrac{1}{2}$ and $m_2+\tfrac{1}{2}$ are non-central critical points for $L(s,\itPi \times \itSigma)$. The existence of $\itPi$ is clear, as it can be constructed as the automorphic induction of some Hecke character over CM-field (cf.\,\cite[Appendix]{RW2004}). 
The existence of $\itSigma$ follows from Lemma \ref{L:auxiliary 1}.
By the result of Raghuram \cite[Theorem 1.1]{Raghuram2009}, for all critical points $m+\tfrac{1}{2}$ for $L(s,\itPi \times \itSigma)$, we have
\begin{align*}
&\sigma \left( \frac{L(m+\tfrac{1}{2},\itPi \times \itSigma)}{p(m,\itPi_\infty \times \itSigma_\infty)\cdot G(\omega_\itSigma)\cdot p(\itPi,(-1)^{m+n}\varepsilon(\itSigma_\infty))\cdot p(\itSigma,\varepsilon(\itSigma_\infty))} \right)\\
& =  \frac{L(m+\tfrac{1}{2},{}^\sigma\!\itPi \times {}^\sigma\!\itSigma)}{p(m,\itPi_\infty \times \itSigma_\infty)\cdot G({}^\sigma\!\omega_\itSigma)\cdot p({}^\sigma\!\itPi,(-1)^{m+n}\varepsilon(\itSigma_\infty))\cdot p({}^\sigma\!\itSigma,\varepsilon(\itSigma_\infty))} ,\quad \sigma \in {\rm Aut}(\C).
\end{align*}
This implies that
\begin{align}\label{E:auxiliary 2 pf 1}
\begin{split}
&\sigma \left( \frac{L(m_1+\tfrac{1}{2},\itPi \times \itSigma)}{L(m_2+\tfrac{1}{2},\itPi \times \itSigma)}\cdot \frac{p(m_2,\itPi_\infty \times \itSigma_\infty)\cdot  p(\itPi,(-1)^{m_2+n}\varepsilon(\itSigma_\infty))}{p(m_1,\itPi_\infty \times \itSigma_\infty)\cdot  p(\itPi,(-1)^{m_1+n}\varepsilon(\itSigma_\infty))} \right)\\
& =  \frac{L(m_1+\tfrac{1}{2},{}^\sigma\!\itPi \times {}^\sigma\!\itSigma)}{L(m_2+\tfrac{1}{2},{}^\sigma\!\itPi \times {}^\sigma\!\itSigma)}\cdot \frac{p(m_2,\itPi_\infty \times \itSigma_\infty)\cdot  p({}^\sigma\!\itPi,(-1)^{m_2+n}\varepsilon(\itSigma_\infty))}{p(m_1,\itPi_\infty \times \itSigma_\infty)\cdot  p({}^\sigma\!\itPi,(-1)^{m_1+n}\varepsilon(\itSigma_\infty))} ,\quad \sigma \in {\rm Aut}(\C).
\end{split}
\end{align}
On the other hand, by the result of Harder and Raghuram \cite[Theorem 7.21]{HR2020}, if either
\[
(m_1,m_2) = (\tfrac{-{\sf w}-{\sf u}+1}{2}, \tfrac{-{\sf w}-{\sf u}-1}{2})
\]
or $m_1$ and $m_2$ are both greater than or both less than $\tfrac{-{\sf w}-{\sf u}}{2}$, then we have
\begin{align}\label{E:auxiliary 2 pf 2}
\begin{split}
&\sigma \left( \frac{L(m_1+\tfrac{1}{2},\itPi \times \itSigma)}{L(m_2+\tfrac{1}{2},\itPi \times \itSigma)}\cdot (\sqrt{-1})^{(m_1-m_2)n(n-1)/2}\cdot\frac{p(\itPi,(-1)^{m_2+n}\varepsilon(\itSigma_\infty))}{p(\itPi,(-1)^{m_1+n}\varepsilon(\itSigma_\infty))} \right)\\
& =  \frac{L(m_1+\tfrac{1}{2},{}^\sigma\!\itPi \times {}^\sigma\!\itSigma)}{L(m_2+\tfrac{1}{2},{}^\sigma\!\itPi \times {}^\sigma\!\itSigma)}\cdot (\sqrt{-1})^{(m_1-m_2)n(n-1)/2}\cdot\frac{p({}^\sigma\!\itPi,(-1)^{m_2+n}\varepsilon(\itSigma_\infty))}{p({}^\sigma\!\itPi,(-1)^{m_1+n}\varepsilon(\itSigma_\infty))} ,\quad \sigma \in {\rm Aut}(\C).
\end{split}
\end{align}
Therefore, in these cases we conclude from (\ref{E:auxiliary 2 pf 1}) and (\ref{E:auxiliary 2 pf 2}) that
\begin{align*}
&\sigma \left(\frac{p(m_1,\itPi_\infty \times \itSigma_\infty)}{p(m_2,\itPi_\infty \times \itSigma_\infty)} \cdot (\sqrt{-1})^{(m_1-m_2)n(n-1)/2} \right) = \frac{p(m_1,\itPi_\infty \times \itSigma_\infty)}{p(m_2,\itPi_\infty \times \itSigma_\infty)} \cdot (\sqrt{-1})^{(m_1-m_2)n(n-1)/2}
\end{align*}
for all $\sigma \in {\rm Aut}(\C)$, which implies (\ref{E:local period}).
In particular, (\ref{E:local period}) holds when $L(s,\itPi \times \itSigma)$ does not admit central critical point, that is, when ${\sf w}$ is odd.
The remain case is when ${\sf w}$ is even, $m_1 > \tfrac{-{\sf w}-{\sf u}}{2}$, and $m_2<\tfrac{-{\sf w}-{\sf u}}{2}$.
In this case, we may assume $\itPi$ and $\itSigma$ are both essentially self-dual by Lemma \ref{L:auxiliary 1}.
More precisely, by the construction of $\itPi$ and $\itSigma$ in $loc.$ $cit.$ as functorial transfers from $\SO_{n+1}(\A)$  and $\Sp_{n-2}(\A)$ respectively, we have $\itPi = \itPi^\vee \otimes |\mbox{ }|_\A^{\sf w}$ and $\itSigma = \itSigma^\vee \otimes |\mbox{ }|_\A^{\sf u}$ with $\omega_\itPi = |\mbox{ }|_\A^{n{\sf w}/2}$ and $\omega_\itSigma = |\mbox{ }|_\A^{(n-1){\sf u}/2}$.
Therefore, $G(\omega_\itPi) \in \Q^\times$, $G(\omega_\itSigma) \in \Q^\times$, and by Theorem \ref{T:ratio under duality} we have
\begin{align}\label{E:auxiliary 2 pf 3}
\begin{split}
&\sigma \left(  \frac{L(m_2+\tfrac{1}{2},\itPi \times \itSigma)}{\varepsilon(0,\itPi_\infty\times\itSigma_\infty,\psi_\infty)\cdot L(-m_2-{\sf w}-{\sf u}+\tfrac{1}{2},\itPi \times \itSigma)} \right)\\
& = \frac{L(m_2+\tfrac{1}{2},{}^\sigma\!\itPi \times {}^\sigma\!\itSigma)}{\varepsilon(0,\itPi_\infty\times\itSigma_\infty,\psi_\infty)\cdot L(-m_2-{\sf w}-{\sf u}+\tfrac{1}{2},{}^\sigma\!\itPi \times {}^\sigma\!\itSigma)},\quad \sigma \in {\rm Aut}(\C).
\end{split}
\end{align}
By (\ref{E:auxiliary 2 pf 1}) for $m_1, m_2$ and (\ref{E:auxiliary 2 pf 2}) for $m_1,-m_2-{\sf w}-{\sf u}$ (both greater than $\tfrac{-{\sf w}-{\sf u}}{2}$), we then deduce from (\ref{E:auxiliary 2 pf 3}) that
\begin{align*}
&\sigma \left(  \frac{p(m_1,\itPi_\infty \times \itSigma_\infty)}{p(m_2,\itPi_\infty \times \itSigma_\infty)} \cdot \varepsilon(0,\itPi_\infty \times \itSigma_\infty,\psi_\infty)\cdot(\sqrt{-1})^{(m_1+m_2+{\sf w}+{\sf u})n(n-1)/2} \right) \\
& = \frac{p(m_1,\itPi_\infty \times \itSigma_\infty)}{p(m_2,\itPi_\infty \times \itSigma_\infty)} \cdot \varepsilon(0,\itPi_\infty \times \itSigma_\infty,\psi_\infty)\cdot(\sqrt{-1})^{(m_1+m_2+{\sf w}+{\sf u})n(n-1)/2},\quad \sigma \in {\rm Aut}(\C),
\end{align*}
which implies 
\[
\frac{p(m_1,\itPi_\infty \times \itSigma_\infty)}{p(m_2,\itPi_\infty \times \itSigma_\infty)} \cdot \varepsilon(0,\itPi_\infty \times \itSigma_\infty,\psi_\infty)\cdot(\sqrt{-1})^{(m_1+m_2+{\sf w}+{\sf u})n(n-1)/2} \in \Q^\times.
\]
Since both ${\sf w}$ and ${\sf u}$ are even, it is easy to see that 
\[
\varepsilon(0,\itPi_\infty \times \itSigma_\infty,\psi_\infty)\cdot(\sqrt{-1})^{(m_1+m_2+{\sf w}+{\sf u})n(n-1)/2} \in (\sqrt{-1})^{(m_1-m_2)n(n-1)/2}\cdot\Q^\times.
\]
This completes the proof.
\end{proof}

\subsection{Main results}\label{SS:main results}

Following is the main result of this paper. We prove a period relation between the Betti--Whittaker periods under duality.

\begin{thm}\label{T:main 1}
Let $\itPi$ be a regular algebraic cuspidal automorphic representation of $\GL_n(\A)$.
Let $\varepsilon \in \{\pm1\}$ if $n$ is even, and $\varepsilon = \varepsilon(\itPi_\infty)$ if $n$ is odd.
Assume the following conditions are satisfied:
\begin{align}\label{E:regularity 1}
\begin{cases}
\min\{\kappa_i\} \geq 3& \mbox{ if $n$ is even},\\
\min\{\kappa_i\} \geq 5& \mbox{ if $n$ is odd},
\end{cases}
\end{align}
and
\begin{align}\label{E:regularity 2}
\begin{cases}
\mbox{$\underline{\kappa}$ is $4$-regular} & \mbox{ if $n$ or ${\sf w}$ is odd},\\
\mbox{$\underline{\kappa}$ is $6$-regular} &\mbox{ if $n$ and ${\sf w}$ are even}.
\end{cases}
\end{align}
Here $(\underline{\kappa};\,{\sf w})$ is the infinity type of $\itPi$.
Then we have 
\begin{align}\label{E:main}
\sigma \left(\frac{p(\itPi,\varepsilon)}{G(\omega_\itPi)^{n-1}\cdot p(\itPi^\vee,\varepsilon)} \right) = \frac{p({}^\sigma\!\itPi,\varepsilon)}{G({}^\sigma\!\omega_\itPi)^{n-1}\cdot p({}^\sigma\!\itPi^\vee,\varepsilon)} ,\quad \sigma \in {\rm Aut}(\C).
\end{align}
\end{thm}

\begin{proof}
We prove assertion (\ref{E:main}) by induction on $n$.
When $n=1$, as we mentioned in Remark \ref{R:n=1}, we have
$p({}^\sigma\!\itPi,\varepsilon) = C$ for all $\sigma \in {\rm Aut}(\C)$ for some $C \in \C^\times$ uniquely determined by $[\itPi_\infty]^\varepsilon$.
By our normalization (\ref{E:normalization}), we also have
$p({}^\sigma\!\itPi^\vee,\varepsilon) = C$ for all $\sigma \in {\rm Aut}(\C)$.
Thus (\ref{E:main}) holds.
Assume (\ref{E:main}) holds for all $1 \leq n' \leq n-1$ subject to the regularity conditions (\ref{E:regularity 1}) and (\ref{E:regularity 2}).
Let $\itPi$ be a regular algebraic cuspidal automorphic representation of $\GL_n(\A)$ satisfying the regularity conditions.
Put $\delta \in \{0,1\}$ such that $\delta \equiv n\,({\rm mod}\,2)$. Let $r = {\lfloor \tfrac{n}{2}\rfloor}$ and $(\underline{\kappa};\,{\sf w})$ be the infinity type of $\itPi$. Define a tuple of integers $\underline{\ell} \in \Z^{r+\delta-1}$ by
\[
\ell_i = \begin{cases}\kappa_i-2 & \mbox{ if $n$ or ${\sf w}$ is odd},\\
\kappa_i-3 & \mbox{ if $n$ and ${\sf w}$ are even},
\end{cases}
\quad 1 \leq i \leq r-1
\]
and $\ell_{r} = \kappa_r-2$ if $n$ is odd. 
Then it is clear that $(\underline{\ell};\,\delta)$ also satisfies conditions (\ref{E:regularity 1}) and (\ref{E:regularity 2}).
Indeed, suppose $n$ is odd, then $n-1$ is even and $\delta=1$. 
We have $\min\{\ell_i\} = \min\{\kappa_i\}-2 \geq 3$. Also $\underline{\ell}$ is $4$-regular since $|\ell_i-\ell_j| = |\kappa_i-\kappa_j|$ by definition.
The other cases can be verified in a similar way.
By \cite[Appendix]{RW2004} for $n$ odd and Lemma \ref{L:auxiliary 1} for $n$ even, there exists a regular algebraic cuspidal automorphic representation $\itSigma$ of $\GL_{n-1}(\A)$ with infinity type $(\underline{\ell};\,\delta)$. 
Consider the Rankin--Selberg $L$-function $L(s,\itPi \times \itSigma)$.
By definition of $\underline{\ell}$, we see that inequality (\ref{E:balanced}) holds, that is, $(\itPi_\infty,\itSigma_\infty)$ is balanced. Together with condition (\ref{E:regularity 1}), we see that 
\[
\min\{|\kappa_i-\ell_j|,|\kappa_i-1|,|\ell_j-1|\} \geq 2.
\]
In particular, $L(s,\itPi \times \itSigma)$ admits non-central critical points by (\ref{E:critical range}).
We fix a non-central critical point $m+\tfrac{1}{2} \in \Z+\tfrac{1}{2}$ for $L(s,\itPi \times \itSigma)$.
By the result of Raghuram \cite{Raghuram2009}, we have
\begin{align}\label{E:main 1 pf 1}
\begin{split}
&\sigma \left( \frac{L(m+\tfrac{1}{2},\itPi \times \itSigma)}{p(m,\itPi_\infty \times \itSigma_\infty)\cdot G(\omega_\itSigma)\cdot p(\itPi,\varepsilon_m)\cdot p(\itSigma,\varepsilon_m')} \right)\\
& =  \frac{L(m+\tfrac{1}{2},{}^\sigma\!\itPi \times {}^\sigma\!\itSigma)}{p(m,\itPi_\infty \times \itSigma_\infty)\cdot G({}^\sigma\!\omega_\itSigma)\cdot p({}^\sigma\!\itPi,\varepsilon_m)\cdot p({}^\sigma\!\itSigma,\varepsilon_m')} ,\quad \sigma \in {\rm Aut}(\C).
\end{split}
\end{align}
Here $\varepsilon_m = \varepsilon(\itPi_\infty)$ if $n$ is odd, $\varepsilon_m' = \varepsilon(\itSigma_\infty)$ if $n$ is even, and $\varepsilon_m\varepsilon_m' = (-1)^{m+n}$.
Note that $-m+\tfrac{1}{2}$ is a critical point for $L(s,\itPi^\vee \times \itSigma^\vee)$. By $loc.$ $cit.$ for $L(s,\itPi^\vee \times \itSigma^\vee)$, we have
\begin{align}\label{E:main 1 pf 2}
\begin{split}
&\sigma \left( \frac{L(-m+\tfrac{1}{2},\itPi^\vee \times \itSigma^\vee)}{p(-m,\itPi_\infty^\vee \times \itSigma_\infty^\vee)\cdot G(\omega_\itSigma)^{-1}\cdot p(\itPi^\vee,\varepsilon_m)\cdot p(\itSigma^\vee,\varepsilon_m')} \right)\\
& =  \frac{L(-m+\tfrac{1}{2},{}^\sigma\!\itPi^\vee \times {}^\sigma\!\itSigma^\vee)}{p(-m,\itPi_\infty^\vee \times \itSigma_\infty^\vee)\cdot G({}^\sigma\!\omega_\itSigma)^{-1}\cdot p({}^\sigma\!\itPi^\vee,\varepsilon_m)\cdot p({}^\sigma\!\itSigma^\vee,\varepsilon_m')} ,\quad \sigma \in {\rm Aut}(\C).
\end{split}
\end{align}
On the other hand, by Theorem \ref{T:ratio under duality}, we have
\begin{align}\label{E:main 1 pf 3}
\begin{split}
&\sigma \left( \frac{L(m+\tfrac{1}{2},\itPi \times \itSigma)}{\varepsilon(0,\itPi_\infty\times\itSigma_\infty,\psi_\infty)\cdot G(\omega_\itPi)^{n-1}\cdot G(\omega_\itSigma)^{n}\cdot L(-m+\tfrac{1}{2},\itPi^\vee \times \itSigma^\vee)} \right)\\
& = \frac{L(m+\tfrac{1}{2},{}^\sigma\!\itPi \times {}^\sigma\!\itSigma)}{\varepsilon(0,\itPi_\infty\times\itSigma_\infty,\psi_\infty)\cdot G({}^\sigma\!\omega_\itPi)^{n-1}\cdot G({}^\sigma\!\omega_\itSigma)^{n}\cdot L(-m+\tfrac{1}{2},{}^\sigma\!\itPi^\vee \times {}^\sigma\!\itSigma^\vee)},\quad \sigma \in {\rm Aut}(\C).
\end{split}
\end{align}
We thus conclude from (\ref{E:main 1 pf 1})-(\ref{E:main 1 pf 3}) that
\begin{align*}
&\sigma \left( \varepsilon(0,\itPi_\infty\times\itSigma_\infty,\psi_\infty)^{-1}\cdot \frac{p(m,\itPi_\infty \times \itSigma_\infty)}{p(-m,\itPi_\infty^\vee \times \itSigma_\infty^\vee)}\right)\\
&\times\sigma \left(\frac{p(\itPi,\varepsilon_m)}{G(\omega_\itPi)^{n-1}\cdot p(\itPi^\vee,\varepsilon_m)}\cdot \frac{p(\itSigma,\varepsilon_m')}{G(\omega_\itSigma)^{n-2}\cdot p(\itSigma^\vee,\varepsilon_m')}  \right)\\
& = \varepsilon(0,\itPi_\infty\times\itSigma_\infty,\psi_\infty)^{-1}\cdot \frac{p(m,\itPi_\infty\times\itSigma_\infty)}{p(-m,\itPi_\infty^\vee \times \itSigma_\infty^\vee)}\\
&\times \frac{p({}^\sigma\!\itPi,\varepsilon_m)}{G({}^\sigma\!\omega_\itPi)^{n-1}\cdot p({}^\sigma\!\itPi^\vee,\varepsilon_m)}\cdot \frac{p(\itSigma,\varepsilon_m')}{G({}^\sigma\!\omega_\itSigma)^{n-2}\cdot p({}^\sigma\!\itSigma^\vee,\varepsilon_m')}, \quad \sigma \in {\rm Aut}(\C).
\end{align*}
Since $\itPi_\infty^\vee = \itPi_\infty \otimes |\mbox{ }|^{-{\sf w}}$ and $\itSigma_\infty^\vee = \itSigma_\infty \otimes |\mbox{ }|^{-\delta}$, by Lemma \ref{E:archi. period comparison} we have
\[
\frac{p(-m,\itPi_\infty^\vee\times\itSigma_\infty^\vee)}{p(-m-\delta-{\sf w},\itPi_\infty\times\itSigma_\infty)} \in \Q^\times.
\]
It then follows from the archimedean period relation Lemma \ref{L:auxiliary 2} that
\[
\frac{p(m,\itPi_\infty\times\itSigma_\infty)}{p(-m,\itPi_\infty^\vee\times\itSigma_\infty^\vee)} \in (\sqrt{-1})^{(2m+{\sf w}+\delta)n(n-1)/2}\cdot \Q^\times.
\]
Also it is easy to see that
\[
\varepsilon(0,\itPi_\infty \times \itSigma_\infty,\psi_\infty) \in (\sqrt{-1})^{({\sf w}+\delta)n(n-1)/2}\cdot\Q^\times.
\]
We refer to \cite[(3.7)]{Knapp1994} for formulas of archimedean $\varepsilon$-factors.
Therefore, we have
\begin{align*}
&\sigma \left(\frac{p(\itPi,\varepsilon_m)}{G(\omega_\itPi)^{n-1}\cdot p(\itPi^\vee,\varepsilon_m)}\cdot \frac{p(\itSigma,\varepsilon_m')}{G(\omega_\itSigma)^{n-2}\cdot p(\itSigma^\vee,\varepsilon_m')}  \right)\\
& = \frac{p({}^\sigma\!\itPi,\varepsilon_m)}{G({}^\sigma\!\omega_\itPi)^{n-1}\cdot p({}^\sigma\!\itPi^\vee,\varepsilon_m)}\cdot \frac{p(\itSigma,\varepsilon_m')}{G({}^\sigma\!\omega_\itSigma)^{n-2}\cdot p({}^\sigma\!\itSigma^\vee,\varepsilon_m')}, \quad \sigma \in {\rm Aut}(\C).
\end{align*}
In particular, the period relation (\ref{E:main}) holds for $p(\itPi,\varepsilon_m)$ if and only if it holds for $p(\itSigma,\varepsilon_m')$.
Since $\itSigma$ also satisfies the regularity conditions, by induction hypothesis, (\ref{E:main}) holds for $p(\itSigma,\varepsilon_m')$.
When $n$ is even, by replacing $\itSigma$ with $\itSigma \otimes \chi$ for some quadratic Hecke character $\chi$ with $\chi_\infty(-1)=-1$, we see that the period relation (\ref{E:main}) also holds for $p(\itPi,-\varepsilon_m')$.
This completes the proof.
\end{proof}

\begin{rmk}
If we replace $[\itPi_\infty]^\varepsilon$ by $C\cdot [\itPi_\infty]^\varepsilon$ for some $C\in\C^\times$, then $[\itPi_\infty^\vee]^{\varepsilon}$ is also replaced by $C\cdot [\itPi_\infty^\vee]^{\varepsilon}$ according to our normalization (\ref{E:normalization}) (note that $\itPi_\infty^\vee = \itPi_\infty\otimes|\mbox{ }|^{-{\sf w}}$).
Therefore, the assertion (\ref{E:main}) is independent of the choice of generators.
\end{rmk}

In the special case when $\itPi$ is of orthogonal type, we have the following corollary of Theorem \ref{T:main 1} on the period relation between Betti--Whittaker periods with opposite signature.

\begin{corollary}\label{C:main}
Let $\itPi$ be a regular algebraic cuspidal automorphic representation of $\GL_{2n}(\A)$. 
Assume $\itPi$ is $\chi$-orthogonal and satisfies the regularity conditions in Theorem \ref{T:main 1}.
Then we have
\begin{align*}
\sigma \left(\frac{p(\itPi,+)}{G(\chi^n\omega_\itPi^{-1})\cdot p(\itPi,-)} \right) = \frac{p({}^\sigma\!\itPi,+)}{G({}^\sigma\!\chi^n\omega_{{}^\sigma\!\itPi}^{-1})\cdot p({}^\sigma\!\itPi,-)} ,\quad \sigma \in {\rm Aut}(\C).
\end{align*}
\end{corollary}

\begin{proof}
By the result of Raghuram and Shahidi \cite{RS2008}, for any algebraic Hecke character $\eta$ and $\varepsilon \in \{\pm1\}$, we have
\begin{align}\label{E:RS}
\sigma \left( \frac{p(\itPi\otimes\eta,\varepsilon)}{G(\eta)^{n(2n-1)}\cdot p(\itPi,\varepsilon\cdot\varepsilon(\eta_\infty))} \right) = \frac{p({}^\sigma\!\itPi\otimes{}^\sigma\!\eta,\varepsilon)}{G({}^\sigma\!\eta)^{n(2n-1)}\cdot p({}^\sigma\!\itPi,\varepsilon\cdot\varepsilon(\eta_\infty))},\quad \sigma \in {\rm Aut}(\C).
\end{align}
Snce $\itPi$ is $\chi$-orthogonal, we have $\itPi = \itPi^\vee\otimes\chi$ and $\varepsilon(\chi_\infty)=-1$ by Lemma \ref{L:auxiliary 3}.
The assertion then follows immediately from Theorem \ref{T:main 1} and (\ref{E:RS}) with $\eta=\chi^{-1}$.
\end{proof}

In \cite{HR2020}, Harder and Raghuram study the rank-one Eisenstein cohomology of $\GL_{2n+n'}$ and prove the algebraicity of successive critical values for $\GL_{2n} \times \GL_{n'}$.
They show that the ratios are algebraic when $n'$ is even, and the algebraicity is expressed in terms of the relative periods when $n'$ is odd.
As a consequence of Corollary \ref{C:main}, the relative periods are essentially trivial when the representations of $\GL_{2n}$ are of orthogonal type.
More precisely, we have the following:
\begin{thm}\label{T:main 2}
Let $\itPi$ and $\itSigma$ be regular algebraic cuspidal automorphic representations of $\GL_{2n}(\A)$ and $\GL_{n'}(\A)$ respectively. 
Assume $\itPi$ is $\chi$-orthogonal and satisfies the regularity condition (\ref{E:regularity 2}) in Theorem \ref{T:main 1}. Let $m_0,m_0+1 \in \Z+\tfrac{n'}{2}$ be critical points for $L(s,\itPi \times \itSigma)$ such that $L(m_0+1,\itPi \times \itSigma) \neq 0$. Then we have
\begin{align*}
&\sigma \left( \frac{L(m_0,\itPi \times \itSigma)}{(\sqrt{-1})^{nn'}\cdot G(\chi^{n}\omega_\itPi^{-1})^{n'}\cdot L(m_0+1,\itPi \times \itSigma)} \right)\\
& =  \frac{L(m_0,{}^\sigma\!\itPi \times {}^\sigma\!\itSigma)}{(\sqrt{-1})^{nn'}\cdot G({}^\sigma\!\chi^{n}\omega_{{}^\sigma\!\itPi}^{-1})
^{n'}\cdot L(m_0+1,{}^\sigma\!\itPi \times {}^\sigma\!\itSigma)} ,\quad \sigma \in {\rm Aut}(\C).
\end{align*}

\end{thm}

\begin{proof}
We only need to consider the case when $n'$ is odd.
By the result of Harder and Raghuram \cite[Theorem 7.21]{HR2020}, we have
\begin{align*}
&\sigma \left(  \frac{L(m_0,\itPi \times \itSigma)}{L(m_0+1,\itPi \times \itSigma)}\cdot (\sqrt{-1})^{n}\cdot \frac{p(\itPi,(-1)^{m_0+1+n+n'/2}\varepsilon(\itSigma_\infty))}{p(\itPi,(-1)^{m_0+n+n'/2}\varepsilon(\itSigma_\infty))}\right)\\
& =  \frac{L(m_0,{}^\sigma\!\itPi \times {}^\sigma\!\itSigma)}{L(m_0+1,{}^\sigma\!\itPi \times {}^\sigma\!\itSigma)}\cdot (\sqrt{-1})^{n}\cdot \frac{p({}^\sigma\!\itPi,(-1)^{m_0+1+n+n'/2}\varepsilon(\itSigma_\infty))}{p({}^\sigma\!\itPi,(-1)^{m_0+n+n'/2}\varepsilon(\itSigma_\infty))},\quad \sigma \in {\rm Aut}(\C).
\end{align*}
Here we interpret the relative periods of $\itPi$ as
\[
(\sqrt{-1})^{n}\cdot \frac{p(\itPi,\varepsilon)}{p(\itPi,-\varepsilon)},\quad \varepsilon \in \{\pm1\}
\]
which follows essentially from definition as mentioned in \cite[\S\,5.2.3]{HR2020}.
The reader is also referred to \cite[Theorem 3.1]{Raghuram2013} where the relative periods are defined without $(\sqrt{-1})^{n}$.
Note that the regularity condition (\ref{E:regularity 1}) is automatically satisfied by the existence of two critical points $m_0$ and $m_0+1$.
The assertion then follows from Corollary \ref{C:main}.
\end{proof}

Let $\F$ be an \'{e}tale real quadratic algebra over $\Q$. We denote by $\infty_1$ and $\infty_2$ the non-zero algebra homomorphisms from $\F$ into $\R$. By the results of Ramakrishnan \cite{Rama2000} and Krishnamurthy \cite{Kris2003}, the Langlands functoriality of the Asai transfer from $\GL_{2}(\A_\F)$ to $\GL_{4}(\A)$ holds.
As a consequence of Theorem \ref{T:main 2}, we have the following:

\begin{corollary}\label{C:main 2}
Let $\pi$ and $\itSigma$ be regular algebraic cuspidal automorphic representations of $\GL_2(\A_\F)$ and $\GL_{n'}(\A)$ respectively. 
Assume the following regularity condition is satisfied:
\[
\min\{\kappa_1,\kappa_2\}\geq \begin{cases}
3 & \mbox{ if $\kappa_1+\kappa_2$ is even},\\
4 & \mbox{ if $\kappa_1+\kappa_2$ is odd},
\end{cases}
\]
where $\kappa_i \geq 2$ is the minimal ${\rm SO}_2(\R)$-weight of $\pi_{\infty_i}$ for $i=1,2$.
Let $m_0,m_0+1 \in \Z+\tfrac{n'-1}{2}$ be critical points for $L(s,{\rm As}(\pi)\times\itSigma)$ such that $L(m_0+1,{\rm As}(\pi)\times \itSigma) \neq 0$. Then we have
\[
\sigma \left( \frac{L(m_0,{\rm As}(\pi)\times\itSigma)}{G(\omega_{\F/\Q})^{n'}\cdot L(m_0+1,{\rm As}(\pi)\times\itSigma)} \right)=  \frac{L(m_0,{\rm As}({}^\sigma\!\pi)\times{}^\sigma\!\itSigma)}{G(\omega_{\F/\Q})^{n'}\cdot L(m_0+1,{\rm As}({}^\sigma\!\pi)\times{}^\sigma\!\itSigma)} ,\quad \sigma \in {\rm Aut}(\C).
\]
Here $\omega_{\F/\Q}$ is the quadratic Hecke character associated to $\F/\Q$ by class field theory.

\end{corollary}

\begin{proof}
There exist ${\sf w}_1,{\sf w}_2 \in \Z$ such that $\kappa_i \equiv {\sf w}_i\,({\rm mod}\,2)$ and $\pi_{\infty_i} = D_{\kappa_i}\otimes|\mbox{ }|^{{\sf w}_i/2}$ for $i=1,2$.
Then ${\rm As}_\infty(\pi_\infty)\otimes |\mbox{ }|^{1/2} \in \Omega(4)$ and its infinity type is equal to 
\[
(\kappa_1+\kappa_2-1,|\kappa_1-\kappa_2|+1;\,{\sf w}_1+{\sf w}_2+1).
\]
Therefore the regularity conditions (\ref{E:regularity 1}) and (\ref{E:regularity 2}) are satisfied by $\itPi = {\rm As}(\pi)\otimes|\mbox{ }|_\A^{1/2}$.
It is well-known that $\pi = \pi^\vee \otimes \omega_\pi$. By Corollary \ref{C:tensor} (resp.\,Corollary \ref{C:Asai}) when $\F=\Q \times \Q$ (resp.\,$\F$ is a field), ${\rm As}(\pi)\otimes|\mbox{ }|_\A^{1/2}$ is regular algebraic and $\omega_\pi\vert_{\A^\times}|\mbox{ }|_\A$-orthogonal. 
Moreover, it is either cuspidal or an isobaric sum of two $\omega_\pi\vert_{\A^\times}|\mbox{ }|_\A$-orthogonal cuspidal automorphic representations.
Also note that the central character of ${\rm As}(\pi)$ is equal to $\omega_\pi\vert_{\A^\times}^2\cdot\omega_{\F/\Q}$.
The assertion then follows from Theorem \ref{T:main 2} with $\itPi = {\rm As}(\pi)\otimes|\mbox{ }|_\A^{1/2}$.
This completes the proof.
\end{proof}

\section{Compatibility with Deligne's conjectures}\label{S:Deligne}

In this section, we show that Theorem \ref{T:main 1} is compatible with Deligne's conjecture \cite[Conjecture 2.8]{Deligne1979} on critical $L$-values for motives.
In \S\,\ref{SS:DR}, we also briefly introduce a paper of Deligne and Raghuram \cite{DR2023} on relation between Deligne's motivic periods.

\subsection{Yoshida's period invariants of motives}
We begin with the period invariant of motives introduced by Yoshida in \cite{Yoshida2001}.
Fix a positive integers $n$. Let $d^\pm$ such that $d^++d^-=n$. We assume $d^+=d^-$ if $n$ is even and $|d^+-d^-|=1$ if $n$ is odd.
For $\underline{a} = (a_1,\cdots,a_n) \in \Z_{\geq0}^n$ and $(k^+,k^-) \in \Z_{\geq0}^2$, we say a polynomial function $f$ on $n$ by $n$ matrices over $\Q$ is of type $\{\underline{a};(k^+,k^-)\}$ if 
\[
f\left( \bp t_1 & 0 & 0 \\ \vdots & \ddots & 0 \\ * & \cdots & t_n \ep g \bp g_+ & 0 \\ 0 & g_-\ep \right) = \prod_{i=1}^nt_i^{a_i}\cdot (\det g_+)^{k^+}\cdot (\det g_-)^{k^-}\cdot f(g)
\]
for all $t_1,\cdots,t_n \in \GL_1$ and $g_\pm \in \GL_{d^\pm}$. 
We call $\{\underline{a};(k^+,k^-)\}$ admissible if it satisfies conditions in \cite[(1.6)-(1.10)]{Yoshida2001}.
In this case, we say $f$ is admissible. Yoshida prove in \cite[Theorem 1]{Yoshida2001} that, up to scalar multiples, there exists a unique non-zero admissible polynomial of a given admissible type. For instance, we have the following admissible types:
\begin{itemize}
\item[(1)] $\{(1,\cdots,1);(1,1)\}$.
\item[(2)] $\{(\overbrace{1,\cdots,1}^{d^+},0\cdots,0);(1,0)\}$.
\item[(3)] $\{(\overbrace{1,\cdots,1}^{d^-},0\cdots,0);(0,1)\}$.
\item[(4)] $\{(\overbrace{2,\cdots,2}^i,\overbrace{1,\cdots,1}^{n-2i},0\cdots,0);(1,1)\}$, $1 \leq i \leq \lfloor \tfrac{n}{2}\rfloor-1$.
\end{itemize}
For type (1), the corresponding admissible polynomial is just the determinant function. For types (2) and (3), the corresponding admissible polynomials are denoted by $f^+$ and $f^-$ respectively. Then $f^+$ (resp.\,$f^-$) is given by the determinant of the upper left (resp.\,upper right) $d^+$ by $d^+$ (resp.\,$d^-$ by $d^-$) submatrix.
For type (4), the corresponding admissible polynomials are denoted by $f_i$ for $1 \leq i \leq \lfloor \tfrac{n}{2}\rfloor-1$. Moreover, Yoshida prove that any admissible polynomial can be expressed uniquely as a monomial in $\det$, $f^\pm$, and $f_i$ for $1 \leq i \leq \lfloor \tfrac{n}{2}\rfloor-1$.
For an admissible polynomial 
\[
f = (\det)^{m_0}\cdot\prod_{i}f_i^{m_i}\cdot (f^+)^{m_+}\cdot (f^-)^{m_-}, 
\]
we define its dual $f^\vee$ by
\[
f^\vee = (\det)^{m_0}\cdot\prod_{i}f_i^{m_i}\cdot (f^+)^{m_-}\cdot (f^-)^{m_+}. 
\]
Let $\varepsilon \in \{\pm1\}$ if $n$ is even and $\varepsilon = d_+-d_-$ if $n$ is odd. We define a specific admissible polynomial $f_{\rm BW}^\varepsilon$ as follows: When $n=1$, let $f_{\rm BW}^\varepsilon=1$. When $n>1$, let
\begin{align}\label{E:motivic BW period}
f_{\rm BW}^\varepsilon = \prod_i f_i\cdot\begin{cases}
f^\varepsilon & \mbox{ if $n$ is even},\\
f^+\cdot f^- & \mbox{ if $n$ is odd}.
\end{cases}
\end{align}
Let $M$ be a regular pure motive over $\Q$ of rank $n$ with coefficients in a number field $\E$. We take $d^\pm$ be the dimension of the $\pm$-eigenspace of the Betti realization of $M$ under the archimedean Frobenius action.
The comparison isomorphism between the Betti and de Rham realizations of $M$ determines a period matrix $X_M$ defined in \cite[(2.7)]{Yoshida2001} which is an $n$ by $n$ matrix with coefficients in $\E\otimes_\Q\C$. It is well-defined as a coset in 
\[
{}^tB_n(\E)\backslash M_{n,n}(\E\otimes_\Q\C)/(\GL_{d^+}(\E)\times\GL_{d^-}(\E)),
\]
where ${}^tB_n$ is the lower triangular Borel subgroup of $\GL_n$. A period invariant of $M$ is the evaluation $f(X_M)$ of some admissible polynomial $f$. As proved by Yoshida \cite[Corollary 2]{Yoshida2001}, period invariants belong to $(\E\otimes_\Q\C)^\times$.
The Deligne's periods $\delta(M)$ and $c^\pm(M)$ are equal to $\det(X_M)$ and $f^\pm(X_M)$. We denote by $c_i(M)=f_i(X_M)$ for $1 \leq i \leq \lfloor \tfrac{n}{2}\rfloor-1$. We call $\delta(M)$, $c^\pm(M)$, and $c_i(M)$ for $1 \leq i \leq \lfloor \tfrac{n}{2}\rfloor-1$ the fundamental periods of $M$.
For the dual motive $M^\vee$ of $M$, we prove the following result on relation between the period invariants under duality and under Tate twist which generalizes \cite[(5.1.7) and (5.1.8)]{Deligne1979}.

\begin{lemma}\label{L:motivic}
Let $f$ be an admissible polynomial of type $\{\underline{a};(k^+,k^-)\}$. We have
\begin{align*}
f^\vee(X_{M^\vee}) = \delta(M)^{-k^+-k^-}\cdot f(X_M).
\end{align*}
For $t \in \Z$, we also have
\[
f(X_{M(t)}) = \left(1\otimes(2\pi\sqrt{-1})\right)^{t(k^+d^++k^-d^-)}\cdot \begin{cases}
f(X_M) & \mbox{ if $t$ is even},\\
f^\vee(X_M) & \mbox{ if $t$ is odd}.
\end{cases}
\]
\end{lemma}

\begin{proof}
By definition, it suffices to prove the assertion for fundamental periods. Indeed, if $f_1$ and $f_2$ are admissible polynomials of types $\{\underline{a}_1;(k_1^+,k_1^-)\}$ and $\{\underline{a}_2;(k_2^+,k_2^-)\}$ respectively, then $f_1f_2$ is admissible of type $\{\underline{a}_1+\underline{a}_2;(k_1^++k_2^+,k_2^-+k_2^-)\}$. Clearly we have $\delta(M^\vee) = \delta(M)^{-1}$. For $c^\pm(M)$, the assertion was proved by Deligne \cite[(5.1.7)]{Deligne1979}.
Fix $1 \leq i \leq r-1$ with $r=\lfloor \tfrac{n}{2}\rfloor$. We consider $c_i(M)$. Let ${\sf w}$ be the weight of $M$ and $H_B(M)$ the Betti--realization of $M$. We have the Hodge decomposition
\begin{align}\label{E:Hodge}
\begin{split}
H_B(M) \otimes_\Q \C =& \bigoplus_{i=1}^r \left(H_B^{(1-\kappa_i+{\sf w})/2,(\kappa_i-1+{\sf w})/2}(M)\oplus H_B^{(\kappa_i-1+{\sf w})/2,(1-\kappa_i+{\sf w})/2}(M)\right)\\
& \oplus \begin{cases}
0 & \mbox{ if $n$ is even},\\
H^{{\sf w}/2,{\sf w}/2}(M) & \mbox{ if $n$ is odd}, 
\end{cases}
\end{split}
\end{align}
for some $\kappa_1 > \cdots >\kappa_r \geq 2$ and $\kappa_i \equiv {\sf w}+1\,({\rm mod}\,2)$.
Let $\kappa_i > \ell > \kappa_{i+1}$. By the construction of Deligne \cite{Deligne1971} and Scholl \cite{Scholl1990}, there exists a regular pure motive $N$ over $\Q$ of rank $2$ with coefficients in $\E'$ such that the Hodge decomposition is given by
\[
H_B(N) \otimes_\Q \C = H_B^{0,\ell-1}(M)\oplus H_B^{\ell-1,0}(M).
\]
By replacing $\E$ and $\E'$ by $\E\cdot\E'$, we may assume $\E = \E'$. 
By Yoshida's computation \cite[Proposition 12]{Yoshida2001} (see also the simplified formula for regular motives by Bhagwat \cite{Bhagwat2014}), the condition $\kappa_i > \ell > \kappa_{i+1}$ implies that
\begin{align}\label{E:motivic pf 1}
\begin{split}
c^\pm(M\otimes N) = c_i(M)\cdot\delta(N)^i \cdot (c^+(N)\cdot c^-(N))^{r-i}\cdot  \begin{cases}
1 & \mbox{ if $n$ is even},\\
c^{\pm\varepsilon}(N) & \mbox{ if $n$ is odd},
\end{cases}
\end{split}
\end{align}
where $\varepsilon = d^+-d^-$ if $n$ is odd.
Similarly, we have
\begin{align}\label{E:motivic pf 2}
\begin{split}
c^\pm(M^\vee\otimes N^\vee) &= c_i(M^\vee)\cdot\delta(N^\vee)^i \cdot (c^+(N^\vee)\cdot c^-(N^\vee))^{r-i}\cdot  \begin{cases}
1 & \mbox{ if $n$ is even},\\
c^{\pm\varepsilon}(N^\vee) & \mbox{ if $n$ is odd},
\end{cases}\\
& = c_i(M^\vee)\cdot \delta(N)^{i-n}\cdot 
\delta(N)^i \cdot (c^+(N)\cdot c^-(N))^{r-i}\cdot  \begin{cases}
1 & \mbox{ if $n$ is even},\\
c^{\pm\varepsilon}(N) & \mbox{ if $n$ is odd},
\end{cases}
\end{split}
\end{align}
On the other hand, we have $(M\otimes N)^\vee = M^\vee \otimes N^\vee$ and $\delta(M\otimes N) = \delta(M)^2\cdot\delta(N)^n$.
Thus 
\[
c^\pm(M^\vee\otimes N^\vee) = \delta(M)^{-2}\cdot\delta(N)^{-n}\cdot c^\pm(M\otimes N).
\]
Comparing with (\ref{E:motivic pf 1}) and (\ref{E:motivic pf 2}), we obtain
\[
c_i(M^\vee) = \delta(M)^{-2}\cdot c_i(M).
\]
The assertion for Tate twist can be proved in a similar way.
This completes the proof.
\end{proof}

We recall the Deligne's periods of tensor product of motives in good position.
Let $M$ and $N$ be regular pure motives over $\Q$ of rank $n$ and $n-1$ respectively with coefficient in $\E$.
We say $M$ and $N$ are in good position if 
\[
\begin{cases}
\kappa_1> \ell_1 > \kappa_2 > \ell_2>\cdots>\kappa_{r-1}>\ell_{r-1}>\kappa_r & \mbox{ if $n=2r$},\\
\kappa_1> \ell_1 > \kappa_2 > \ell_2>\cdots>\kappa_{r-1}>\ell_{r-1}>\kappa_r>\ell_r& \mbox{ if $n=2r+1$},
\end{cases}
\]
where the integers $\kappa_i$'s and $\ell_j$'s are determined by the Hodge decompositions of $M$ and $N$ as in (\ref{E:Hodge}).
The following lemma is a direct consequence of Yoshida's computation \cite[Proposition 12]{Yoshida2001}.
\begin{lemma}\label{L:balanced}
Assume $M$ and $N$ are in good position. We have
\[
c^\pm(M\otimes N) = \delta(N)\cdot f_{\rm BW}^\varepsilon(X_M)\cdot f_{\rm BW}^{\varepsilon'}(X_N),
\]
where $\varepsilon = d_M^+-d_M^-$ if $n$ is odd, $\varepsilon' = d_N^+-d_N^-$ if $n$ is even, and $\varepsilon\varepsilon'=\pm1$.
\end{lemma}

\subsection{Compatibility with Deligne's conjecture}

Let $\itPi$ be a regular algebraic cuspidal automorphic representation of $\GL_n(\A)$ with infinity type $(\underline{\kappa};\,{\sf w})$. In \cite[Conjecture 4.5]{Clozel1990}, Clozel proposed a conjecture on the existence of a motive $M_\itPi$ over $\Q$ of rank $n$ with coefficients in some number field $\E$ containing $\Q(\itPi)$ satisfying the following properties:
\begin{itemize}
\item[(1)] $M_\itPi$ has weight $-{\sf w}-n+1$.
\item[(2)] $L(M_\itPi,s) = \left( L^{(\infty)}(s+\tfrac{n-1}{2},{}^\sigma\!\itPi)\right)_{\sigma:\E\rightarrow \C}$.
\item[(3)] $L_\infty(M_\itPi,s) = L(s+\tfrac{n-1}{2},\itPi_\infty)$.
\end{itemize}
Here $L^{(\infty)}$ denotes $L$-function obtained by excluding the archimedean factor. Note that condition (3) implies that $M_\itPi$ is pure and we have the Hodge decomposition
\begin{align*}
H_B(M_\itPi) \otimes_\Q \C =& \bigoplus_{i=1}^r \left(H_B^{(-\kappa_i-{\sf w}-n+2)/2,(\kappa_i-{\sf w}-n)/2}(M_\itPi)\oplus H_B^{(\kappa_i-{\sf w}-n)/2,(-\kappa_i-{\sf w}-n+2)/2}(M_\itPi)\right)\\
& \oplus \begin{cases}
0 & \mbox{ if $n$ is even},\\
H^{(-{\sf w}-n+1)/2,(-{\sf w}-n+1)/2}(M_\itPi) & \mbox{ if $n$ is odd}.
\end{cases}
\end{align*}
Also we have $\varepsilon(\itPi_\infty) = d_{M_\itPi}^+-d_{M_\itPi}^-$.
In the following proposition, under Clozel's and Deligne's conjectures, we show that the automorphic period relation (\ref{E:main}) is a consequence of the motivic period relation Lemma \ref{L:motivic}. 
\begin{prop}
The period relation (\ref{E:main}) holds under the following assumptions:
\begin{itemize}
\item The validity of Clozel's and Deligne's conjectures \cite[Conjecture 4.5]{Clozel1990} and \cite[Conjecture 2.8]{Deligne1979}.
\item The regularity conditions (\ref{E:regularity 1}) and (\ref{E:regularity 2}) hold.
\end{itemize} 
\end{prop}

\begin{proof}
We connect the Betti--Whittaker periods with the motivic period invariants through the Rankin--Selberg $L$-functions for $\GL_n \times \GL_{n-1}$.
Let $\itSigma$ be a regular algebraic cuspidal automorphic representation of $\GL_{n-1}(\A)$ such that $(\itPi_\infty,\itSigma_\infty)$ is balanced and $L(s,\itPi\times\itSigma)$ admits a non-central critical point $m+\tfrac{1}{2}$.
The existence of $\itSigma$ follows from \cite[Appendix]{RW2004} for $n$ odd and Lemma \ref{L:auxiliary 1} for $n$ even, together with the regularity conditions.
Let $M_\itPi$ and $M_{\itSigma}$ be the associate regular pure motives over $\Q$ with coefficients in $\E$.
In the following equalities, the subscript $\sigma$ runs through complex embeddings of $\E$.
Note that 
\[
L(M_\itPi \otimes M_\itSigma,s) = \left(L^{(\infty)}(s+n+\tfrac{1}{2},{}^\sigma\!\itPi \times {}^\sigma\!\itSigma)\right)_\sigma.
\]
By the result of Raghuram \cite{Raghuram2009}, we have
\begin{align*}
L(M_\itPi \otimes M_\itSigma,m-n) \in & \left(1\otimes \frac{p(m,\itPi_\infty\times\itSigma_\infty)}{L(m+\tfrac{1}{2},\itPi_\infty\times\itSigma_\infty)}\right)\cdot \left( G({}^\sigma\!\omega_\itSigma)\cdot p({}^\sigma\!\itPi,\varepsilon)\cdot p({}^\sigma\!\itSigma,\varepsilon')\right)_\sigma\cdot \E^\times,
\end{align*}
where $\varepsilon = \varepsilon(\itPi_\infty)$ if $n$ is odd, $\varepsilon' = \varepsilon(\itSigma_\infty)$ if $n$ is even, and $\varepsilon\varepsilon'=\pm=(-1)^{m+n}$.
By Deligne's conjecture, we have
\[
\frac{L(M_\itPi \otimes M_\itSigma,m-n)}{c^\pm(M_\itPi \otimes M_\itSigma)} \in \left( 1 \otimes (2\pi\sqrt{-1})\right)^{(m-n)n(n-1)/2}\cdot\E^\times.
\]
It is easy to see that $(\itPi_\infty,\itSigma_\infty)$ is balanced if and only if $M_\itPi$ and $M_\itSigma$ are in good position. 
Also note that
\begin{align*}
\delta(M_\itPi) &= \left( 1 \otimes (2\pi\sqrt{-1})\right)^{n{\sf w}/2+n(n-1)/2}\cdot(G({}^\sigma\!\omega_\itPi))_{\sigma},\\
\delta(M_\itSigma) &= \left( 1 \otimes (2\pi\sqrt{-1})\right)^{(n-1){\sf u}/2+(n-1)(n-2)/2}\cdot(G({}^\sigma\!\omega_\itSigma))_{\sigma}.
\end{align*}
Here $|\omega_\itPi| = |\mbox{ }|_\A^{n{\sf w}/2}$ and $|\omega_\itSigma| = |\mbox{ }|_\A^{(n-1){\sf u}/2}$.
We thus deduce from Lemma \ref{L:balanced} that 
\begin{align}\label{E:compatibility pf 1}
\begin{split}
\frac{\left( p({}^\sigma\!\itPi,\varepsilon)\cdot p({}^\sigma\!\itSigma,\varepsilon') \right)_\sigma}{f_{\rm BW}^\varepsilon(X_{M_\itPi})\cdot f_{\rm BW}^{\varepsilon'}(X_{M_\itSigma})} \in &  \left( 1\otimes(2\pi\sqrt{-1})\right)^{(m-n)n(n-1)/2+(n-1){\sf u}/2+(n-1)(n-2)/2}\\
&\times\left( 1\otimes \frac{L(m+\tfrac{1}{2},\itPi_\infty\times\itSigma_\infty)}{p(m,\itPi_\infty\times\itSigma_\infty)} \right)\cdot \E^\times.
\end{split}
\end{align}
By Lemmas \ref{E:archi. period comparison} and \ref{L:auxiliary 2}, we have
\[
\frac{p(m,\itPi_\infty \times \itSigma_\infty)}{p(-m,\itPi_\infty^\vee \times \itSigma_\infty^\vee)} \in (\sqrt{-1})^{(2m+{\sf w}+{\sf u})n(n-1)/2}\cdot \Q^\times.
\]
Also 
\[
\frac{L(m+\tfrac{1}{2},\itPi_\infty\times\itSigma_\infty)}{L(-m+\tfrac{1}{2},\itPi_\infty^\vee\times\itSigma_\infty^\vee)} \in \pi^{-(2m+{\sf w}+{\sf u})n(n-1)/2}\cdot\Q^\times.
\]
Since (\ref{E:compatibility pf 1}) also holds with $\itPi$, $\itSigma$, $m$, ${\sf u}$ replaced by $\itPi^\vee$, $\itSigma^\vee$, $-m$, $-{\sf u}$, we conclude that
\begin{align*}
&\left(\frac{p({}^\sigma\!\itPi,\varepsilon)}{p({}^\sigma\!\itPi^\vee,\varepsilon)}\right)_\sigma  \cdot  \frac{f_{\rm BW}^\varepsilon(X_{M_{\itPi^\vee}})}{f_{\rm BW}^\varepsilon(X_{M_\itPi})} \cdot \left( 1\otimes(2\pi\sqrt{-1})\right)^{{\sf w}n(n-1)/2}\\
&\times \left(\frac{p({}^\sigma\!\itSigma,\varepsilon)}{p({}^\sigma\!\itSigma^\vee,\varepsilon)}\right)_\sigma  \cdot  \frac{f_{\rm BW}^{\varepsilon'}(X_{M_{\itSigma^\vee}})}{f_{\rm BW}^{\varepsilon'}(X_{M_\itSigma})} \cdot \left( 1\otimes(2\pi\sqrt{-1})\right)^{{\sf u}(n-1)(n-2)/2}
\end{align*}
belongs to $\E^\times$. Note that $\itSigma$ also satisfies the regularity conditions. Therefore, an inductive argument shows that 
\[
\left(\frac{p({}^\sigma\!\itPi,\varepsilon)}{p({}^\sigma\!\itPi^\vee,\varepsilon)}\right)_\sigma  \cdot  \frac{f_{\rm BW}^\varepsilon(X_{M_{\itPi^\vee}})}{f_{\rm BW}^\varepsilon(X_{M_\itPi})} \cdot \left( 1\otimes(2\pi\sqrt{-1})\right)^{{\sf w}n(n-1)/2} \in \E^\times.
\]
This implies that (\ref{E:main}) holds if and only if
\[
f_{\rm BW}^\varepsilon(X_{M_{\itPi^\vee}}) = \left( 1\otimes(2\pi\sqrt{-1})\right)^{n(n-1)^2/2}\cdot\delta(M_\itPi)^{-n+1}\cdot f_{\rm BW}^\varepsilon(X_{M_\itPi}).
\]
Since $M_{\itPi^\vee} = (M_\itPi)^\vee(n-1)$, the above period relation follows immediately from \ref{L:motivic}.
This completes the proof.
\end{proof}

\subsection{Period relations of Deligne and Raghuram}\label{SS:DR}

In the paper \cite{DR2023} of Deligne and Raghuram, period relation between $c^+(M)$ and $c^-(M)$ for a pure motive $M$ over $\Q$ is studied.
A criterion is given in \cite[Theorem 2.6.1]{DR2023} under which $c^+(M)$ and $c^-(M)$ are equal up to algebraic multiples. Explicit examples are represented in \cite[\S\,3]{DR2023} including Asai motives and orthogonal or symplectic motives. 
For example, if $M=\otimes_{\F/\Q}M_\pi$ is the (conjectural) Asai motive associated to a regular algebraic cuspidal automorphic representation $\pi$ of $\GL_2(\A_\F)$ for some real quadratic field $\F=\Q(\sqrt{D})$, then it is proved in \cite[Theorem 3.6.2]{DR2023} (special case when $K/F = \F/\Q$) that 
\[
\frac{c^+(M)}{c^-(M)} \in \sqrt{D}\cdot\E^\times.
\]
This is compatible with Corollary \ref{C:main 2} since $G(\omega_{\F/\Q}) \in \sqrt{D}\cdot\Q^\times$.
We also refer to \cite[\S\,3.7.6]{DR2023} for the compatibility between the result \cite{BR2020} and \cite[Theorem 3.7.2]{DR2023} for orthogonal motives.
For general $\chi$-orthogonal regular algebraic cuspidal automorphic representation $\itPi$ of $\GL_{2n}(\A)$, it would be interesting to carry out the period relation between $c^+(M_\itPi)$ and $c^-(M_\itPi)$ by following the arguments in $loc.$ $cit.$. By Theorem \ref{T:main 2}, if Deligne's conjecture holds for $L(M_\itPi,s)$, then we have
\[
\frac{c^+(M_\itPi)}{c^-(M_\itPi)} \in (G({}^\sigma\!\chi{}^\sigma\!\omega_\itPi^{-1}))_\sigma\cdot \E^\times.
\]

\subsection*{Funding}

This work was supported by the Japan Society for the Promotion of Science [JSPS KAKENHI Grant Number 22F22316]; and the National Science and Technology Council of R.O.C. [113-2115-M-007-002-MY3].


\subsection*{Acknowledgement}

The author would like to thank Yao Cheng, Atsushi Ichino, Chi-Heng Lo, and Masao Oi for helpful discussions, and 
Anantharam Raghuram for valuable comments and suggestions on the first draft of this paper.
Thanks are also due to the referee for his/her careful reading and valuable comments and suggestions, especially for pointed out the beautiful paper \cite{DR2023}.




\end{document}